\documentclass{elsart}
\journal{arXiv.org}

\usepackage{amsfonts,amssymb,amsmath,amscd,xy}
\newcommand{\CC}{\mathbb{C}}
\newcommand{\NN}{\mathbb{N}}
\newcommand{\N}{\mathbb{N}}

\newcommand{\Der}{\mathrm{Der}}
\newcommand{\DDer}{\mathbb{D}\mathrm{er}}

\newcommand{\lck}[1]{*+[F]{#1}}
\newcommand{\rep}{\mathsf{rep}}
\newcommand{\iss}{\mathsf{iss}}
\newcommand{\GL}{\mathsf{GL}}

\newcommand{\db}[1]{\{\!\!\{ #1 \}\!\!\}}
\newcommand{\Cxy}{\CC\langle x,y\rangle}

\newcommand{\Neck}{\mathfrak{n}}
\newcommand{\gl}{\mathfrak{gl}}

\newcommand{\cO}{\mathcal{O}}
\newcommand{\DD}{T_A} 

\newtheorem{theorem}{Theorem}
\newtheorem{prop*}{Proposition}
\newtheorem{theorem*}{Theorem}
\newtheorem{lemma}{Lemma}

\newcommand{\prf}{\noindent\textbf{Proof.}~}
\newcommand{\eop}{\hfill $\square$\\}

\newenvironment{proof}{\prf }{\eop}

\theoremstyle{definition}

\newtheorem{remark}{Remark}

\xyoption{all}
\xyoption{poly}
\allowdisplaybreaks

\begin{document}

\begin{frontmatter}
\title{Double Poisson Cohomology of Path Algebras of Quivers}
\author{Anne Pichereau\thanksref{marie}}
\address{Anne Pichereau\\ Department of Mathematics and Computer Science\\ University of Antwerp\\ B-2020 Antwerpen (Belgium)}
\ead{anne.pichereau@ua.ac.be}
\author{Geert Van de Weyer\thanksref{fwo}
}
\address{Geert Van de Weyer\\ Department of Mathematics and Computer Science\\ University of Antwerp\\ B-2020 Antwerpen (Belgium)}
\ead{geert.vandeweyer@ua.ac.be}
\thanks[marie]{
The first author was supported by a grant in the Marie Curie Research Training Network MRTN-CT 2003-505078.}
\thanks[fwo]{
The second author is postdoctoral fellow of the Fund for Scientific Research \--- Flanders (F.W.O.-Vlaanderen)(Belgium).}

\begin{abstract}
In this note, we give a description of the graded Lie algebra of
double derivations of a path algebra as a graded version of the
necklace Lie algebra equipped with the Kontsevich
bracket. Furthermore, we formally introduce the notion of double
Poisson-Lichnerowicz cohomology for double Poisson algebras, and give
some elementary properties. We introduce the notion of a linear double
Poisson tensor on a quiver and show that it induces the structure of a
finite dimensional algebra on the vector spaces $V_v$ generated by the
loops in the vertex $v$. We show that the Hochschild cohomology of the
associative algebra can be recovered from the double Poisson
cohomology. Then, we use the description of the graded necklace Lie
algebra to determine the low-dimensional double Poisson-Lichnerowicz
cohomology groups for three types of (linear and non-linear) double
Poisson brackets on the free algebra $\Cxy$. This allows
us to develop some useful techniques for the computation of the double
Poisson-Lichnerowicz cohomology.
\end{abstract}

\end{frontmatter}

\section{Introduction}\label{intro}
Throughout this paper we will work over an algebraically closed field of characteristic $0$ which we denote by $\CC$. Unadorned tensor products will be over $\CC$. We will use Sweedler notation to write down elements in the tensor product $A\otimes A$ for $A$ an algebra over $\CC$.

A double Poisson algebra $A$ is an associative unital algebra equipped with a linear map
$$\db{-,-}:A\otimes A\rightarrow A\otimes A$$
that is a derivation in its second argument for the outer $A$-bimodule structure on $A\otimes A$, where the outer action of $A$ on $A\otimes A$ is defined as $a.(a'\otimes a'').b := (aa')\otimes (a''b)$. Furthermore, we must have that $\db{a,b} = -\db{b,a}^o$ and that the double Jacobi identity holds for all $a,b,c\in A$:
\begin{align*}
&\db{a,\db{b,c}'}\otimes\db{b,c}'' + \db{c,a}''\otimes \db{b,\db{c,a}'}\\ & + \db{c,\db{a,b}'}''\otimes\db{a,b}''\otimes\db{c,\db{a,b}'}'
= 0,
\end{align*}
where we used Sweedler notation, that is $\db{x,y} = \sum
\db{x,y}'\otimes\db{x,y}''$ for all $x,y\in A$. Such a map is called a
\emph{double Poisson bracket} on $A$. 

Double Poisson algebras were introduced in \cite{MichelDPA} as a
generalization of classical Poisson geometry to the setting of
noncommutative geometry. More specifically, a double Poisson bracket
on an algebra $A$ induces a Poisson structure on all finite
dimensional representation spaces $\rep_n(A)$ of this algebra. Recall
that the coordinate ring $\CC[\rep_n(A)]$ is generated as a
commutative algebra by the generators $a_{ij}$ for $a\in A$ and $1\leq
i,j \leq n$, subject to the relations $\sum_j a_{ij}b_{jk} = (ab)_{ik}$. For each $n$, the Poisson bracket on the coordinate ring $\CC[\rep_n(A)]$ of the variety of $n$-dimensional representations of $A$ is defined as $\{a_{ij},b_{k\ell}\} := \db{a,b}_{kj}'\db{a,b}_{i\ell}''$. This bracket restricts to a Poisson bracket on $\CC[\rep_n(A)]^{\GL_n}$, the coordinate ring of the quotient variety $\iss_n(A)$ under the action of the natural symmetry group $\GL_n$ of $\rep_n(A)$.

In case the algebra is formally smooth (i.e. quasi-free in the sense of \cite{CQ}), double Poisson brackets are completely determined by double Poisson tensors, that is, degree two elements in the tensor algebra $T_A\Der(A,A\otimes A)$. For example, the classical double Poisson bracket on the double $\overline{Q}$ of a quiver $Q$ is the bracket corresponding to the double Poisson tensor 
$$P_{sym} = \sum_{a\in Q}\frac{\partial}{\partial a}\frac{\partial}{\partial a^*}$$
and its Poisson bracket corresponds to the symplectic form on the
representation space of the double of a quiver used in the study of
(deformed) preprojective algebras (see \cite{Bill} and references
therein for further details on deformed preprojective algebras).

We will denote by $\DDer(A)$ the space $\Der(A,A\otimes A)$ of all
derivations of $A$, with value in $A\otimes A$, for the outer
$A$-bimodule structure on $A\otimes A$. This space $\DDer(A)$ becomes a
$A$-bimodule, by using the inner $A$-bimodule structure on $A\otimes A$: if
$\delta\in\DDer(A)$ and $a,b,c\in A$, then $(a\delta
b)(c)=\delta(c)'b\otimes a\delta(c)''$.

As in the classical case, it is possible to define Poisson cohomology for a double Poisson bracket. This was briefly mentioned in \cite{DPSSA} and will be formalized and illustrated in this note. More specifically, in Section \ref{prelims}, we will recall and formalize the definition of the double Poisson cohomology from \cite{DPSSA}. We will then give, in Section \ref{NCMultiQ}, an explicit formulation of the Gerstenhaber algebra of poly-vectorfields and its noncommutative Schouten bracket for the path algebra of a quiver in terms of its graded necklace Lie algebra equipped with a graded version of the Kontsevich bracket. This description will first of all be used to define and classify linear double Poisson structures on path algebras and quivers in Section \ref{DoubleLiePoisson}. On the free algebra in $n$ variables, treated in Section~\ref{example}, this classification becomes
\begin{prop*}[Prop. \ref{dpass1to1}, Section \ref{example}]
There is a one-to-one correspondence between linear double Poisson brackets on $\CC\langle x_1,\dots,x_n\rangle$ and associative algebra structures on $V = \CC x_1\oplus\dots\oplus\CC x_n$. Explicitly, consider the associative algebra structure on $V$ determined by 
$$x_ix_j := \sum_{i,j,k=1}^n c_{ij}^k x_k,$$
where $c_{ij}^k\in\CC$, for all $1\leq i,j,k\leq n$,
then the corresponding double Poisson bracket is given by
$$\db{x_i,x_j} = \sum_{k=1}^n (c_{ij}^k\, x_k\otimes 1 - c_{ji}^k\, 1 \otimes x_k),$$
\end{prop*}
which corresponds to the Poisson tensor:
$$P = \sum_{i,j,k=1}^n c_{ij}^kx_k\frac{\partial}{\partial x_i}\frac{\partial}{\partial x_j}.$$
Next we show there is a connection between the Hochschild cohomology of finite dimensional algebras and the double Poisson cohomology of linear double Poisson structures. We obtain
\begin{theorem*}[Thm. \ref{hochschild}, Section \ref{example}]
Let $A = \CC x_1 \oplus \dots \oplus \CC x_n$ be an $n$-dimensional vector space and let 
$$P = \sum_{i,j,k=1}^n c_{ij}^kx_k\frac{\partial}{\partial x_i}\frac{\partial}{\partial x_j}$$
be a linear double Poisson structure on $T_\CC A = \CC\langle
x_1,\dots, x_n\rangle$. Consider $A$ as an algebra through the product
induced by the structure constants of $P$ (the $c_{ij}^k\in\CC$) and let $HH^\bullet(A)$ denote the Hochschild cohomology of this algebra, then
$$(H_P^\bullet(T_\CC A))_1 \cong HH^\bullet(A).$$
Here the grading on $(H_P^\bullet(T_\CC A))$ is induced by the grading on $\left (\frac{T_{T_\CC A}}{[T_{T_\CC A},T_{T_\CC A}]}\right )_{i}$, which is defined through $\deg(x_i)  =1$.
\end{theorem*}

From the appendix in \cite{MichelDPA} we know that the double Poisson
cohomology of a double Poisson bracket corresponding to a
bi-symplectic form (as defined in \cite{CBEG}) is equal to the
noncommutative de Rham cohomology computed in \cite{RafLieven}, which
is a translation of a similar result in classical Poisson geometry. In
general, little is known about the classical Poisson cohomology and it
is known to be hard to compute. In Section \ref{H0H1}, we will
compute, using the description of the algebra of poly-vectorfields in
Section \ref{NCMultiQ}, the low-dimensional double Poisson cohomology
groups for the free algebra $\CC\langle x,y\rangle$, equipped with
three different types of non-symplectic double Poisson brackets. This
will in particular allow us to develop some tools (including a noncommutative Euler
formula, Proposition \ref{nceuler}) and techniques, that seem to be
useful for the determination of the double Poisson cohomology.

\section{Double Poisson Cohomology}\label{prelims}
In \cite{Lich}, Lichnerowicz observed that $d_\pi = \{\pi,-\}$ with
$\pi$ a Poisson tensor for a Poisson manifold $M$ ($\{-,-\}$ is the Schouten-Nijenhuis bracket) is a square zero differential of degree $+1$, which yields a complex
$$0\stackrel{d_\pi}{\rightarrow} \mathcal{O}(M) \stackrel{d_\pi}{\rightarrow} \Der(\mathcal{O}(M)) \stackrel{d_\pi}{\rightarrow} \wedge^2 \Der(\mathcal{O}(M))\stackrel{d_\pi}{\rightarrow} \dots,$$
the homology of which is called the \emph{Poisson-Lichnerowicz cohomology}. In this section, we show there is an analogous cohomology on $T_A\DDer(A)$ that descends to the classical Poisson-Lichnerowicz cohomology on the quotient spaces of the representation spaces of the algebra.

In \cite[\S 4]{MichelDPA}, the notion of a differentiable double Poisson algebra was introduced. For an algebra $A$, the noncommutative analogue of the classical graded Lie algebra $(\bigwedge_{\cO(M)}\Der(\cO(M)),\{-,-\})$ of polyvector fields on a manifold $M$, where $\{-,-\}$ is the Schouten-Nijenhuis bracket, is the graded Lie algebra
$$T_A\DDer(A)/[T_A\DDer(A),T_A\DDer(A)][1]$$
with graded Lie bracket $\{-,-\} := \mu_A\circ\db{-,-}$ where $\mu_A$ is the multiplication map on $A$ and $\db{-,-}$ is the \emph{double Schouten bracket} defined in \cite[\S 3.2]{MichelDPA}. The classical notion of a Poisson tensor in this new setting becomes
\begin{prop}[{\cite[\S 4.4]{MichelDPA}}]
Let $P\in (T_A\DDer(A))_2$ such that $\{P,P\} = 0$, then $P$ determines a double Poisson bracket on $A$. We call such elements \emph{double Poisson tensors}.
\end{prop}
In case $A$ is formally smooth (for example if $A$ is a path
algebra of a quiver), there is a one-to-one correspondence between
double Poisson tensors on $A$ and double Poisson brackets on $A$.
For a double Poisson tensor $P=\delta\Delta$, the corresponding
double Poisson bracket is, for $a,b\in A$, determined by
$$\db{a,b}_P = \delta(a)'\Delta\delta(a)''(b)-\Delta(a)'\delta\Delta(a)''(b).$$
In order to obtain the noncommutative analogue of the classical Poisson-Lichnerowicz cohomology, we observe that $T_A/[T_A,T_A][1]$ is a graded Lie algebra, so it is a well-known fact that if $P$ satisfies $\{P,P\} = 0$ the map
$$d_P := \{P,-\} : T_A/[T_A,T_A][1]\rightarrow T_A/[T_A,T_A][1]$$
is a square zero differential of degree $+1$. This leads to
\begin{defn}
Let $A$ be a differentiable double Poisson algebra with double Poisson tensor $P$, then the homology $H_P^\bullet(A)$ of the complex
$$0\stackrel{d_P}{\rightarrow}T_A/[T_A,T_A][1]_0\stackrel{d_P}{\rightarrow} T_A/[T_A,T_A][1]_1\stackrel{d_P}{\rightarrow}T_A/[T_A,T_A][1]_2\stackrel{d_P}{\rightarrow}\dots$$
is called the \emph{double Poisson-Lichnerowicz cohomology of $A$}.
\end{defn}

Analogous to the classical interpretation of the first Poisson cohomology groups, we have the following interpretation of the double Poisson cohomology groups:
\begin{eqnarray*}
H^0_P(A) & = & \{\mathrm{double~Casimir~functions}\}\\
         &:= &\{a\in A \mid a\mod [A,A]\in Z(A/[A,A])\}\\
H^1_P(A) & = & \{\mathrm{double~Poisson~vector~fields}\}/\{\mathrm{double~Hamiltonian~vector~fields}\},
\end{eqnarray*}
where in analogy to the classical definitions, a double Poisson vector field is a degree $1$ element $\delta\in T_A/[T_A,T_A]$ satisfying $\{P,\delta\} = 0$ and a double Hamiltonian vector field is a degree $1$ element of the form $\{P,f\}$ with $f\in A/[A,A]$.
Indeed, let us illustrate the first claim. We have for $a\in A$ that
$$\{\delta\Delta,a\} = +\Delta(a)'\delta\Delta(a)'' - \delta(a)'\Delta\delta(a)''$$
whence for any $P\in (T_A\DDer(A))_2$ we get
$$\{P,a\}(b) = - \db{a,b}_{P}$$
so if $P$ is a double Poisson tensor and this expression is zero modulo commutators then $a \mod [A,A]$ is indeed a central element of the Lie algebra $(A/[A,A],\{-,-\}_P)$.

Let $A$ be an associative algebra with unit. From \cite[\S 7]{MichelDPA} we know that the Poisson bracket on $\rep_{n}(A)$ and $\iss_n(A)$ induced by a double Poisson tensor $P$ corresponds to the Poisson tensor $tr(P)$. We furthermore know that the map $tr:\DD/[\DD,\DD][1]\rightarrow\bigwedge\Der(\mathcal{O}(\rep_n(A))$ is a morphism of graded Lie algebras, so we have a morphism of complexes
$$\xymatrix{
0 \ar[r] & (\DD/[\DD,\DD])_0\ar[d]^{tr}\ar[r]^{d_P} &(\DD/[\DD,\DD])_1\ar[d]^{tr}\ar[r]^{d_P} &(\DD/[\DD,\DD])_2\ar[d]^{tr}\ar[r] & \dots \\
0 \ar[r] & \mathcal{O}(\rep_n(A)) \ar[r]^{d_{tr(P)}} & \Der(\mathcal{O}(\rep_n(A))) \ar[r]^{d_{tr(P)}} & \wedge^2\Der(\mathcal{O}(\rep_n(A))) \ar[r] & \dots
}$$
which restricts to a morphism of complexes
$$\xymatrix{
0 \ar[r] & (\DD/[\DD,\DD])_0\ar[d]^{tr}\ar[r]^{d_P} &(\DD/[\DD,\DD])_1\ar[d]^{tr}\ar[r]^{d_P} &(\DD/[\DD,\DD])_2\ar[d]^{tr}\ar[r] & \dots \\
0 \ar[r] & \mathcal{O}(\iss_n(A)) \ar[r]^{d_{tr(P)}} & \Der(\mathcal{O}(\iss_n(A))) \ar[r]^{d_{tr(P)}} & \wedge^2\Der(\mathcal{O}(\iss_n(A))) \ar[r] & \dots
}$$
So there is a map from the double Poisson-Lichnerowicz cohomology to the classical Poisson-Lichnerowicz cohomology on $\rep_n(A)$ and $\iss_n(A)$.

\begin{remark}
It is a well-known fact that in classical Poisson cohomology, because
of the biderivation property satisfied by the Poisson bracket, the
Casimir elements form an algebra. The higher order cohomology groups
all are modules over this algebra. Note that in case of double Casimir
elements, this no longer is the case, as for two double Casimir elements $f$ and $g$ of a double Poisson tensor $P$, we have
$$\db{P,fg} = f\db{P,g} + \db{P,f}g \mathrm{~whence~} 
\{P,fg\} \in f[A,A] + [A,A]g \not\subseteq [A,A].$$
\end{remark}

It is a natural and interesting question to ask whether the map from the double Poisson cohomology to the classical double Poisson cohomology is onto or not. 
For finite-dimensional semi-simple algebras it is onto (see
\cite{DPSSA}), but for some of the double Poisson brackets considered
in Section \ref{H0H1}, the map is not onto. 

\section{NC Multivector Fields and the graded Necklace Lie Algebra}\label{NCMultiQ}
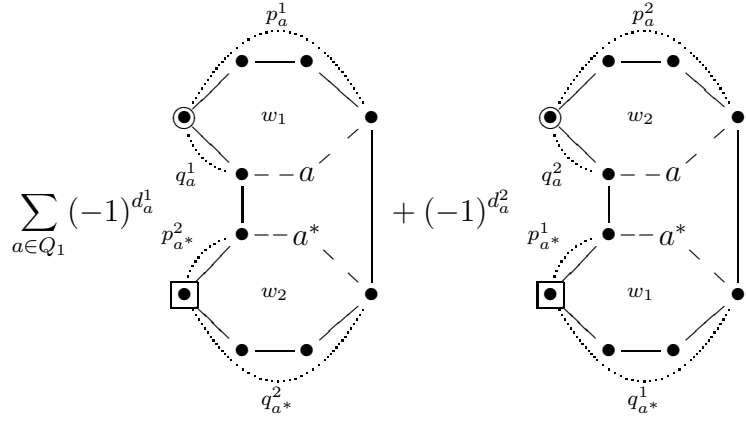
\begin{figure}\label{gradedKbracket}
$$
\sum_{a\in Q_1} (-1)^{d_a^1}
\raisebox{5pc}{\xymatrix@R=1pc@C=1pc{
& \bullet \ar@{-}[r] & \bullet \\
*+[o][F-]{\bullet}\ar@{-}[ur]\ar@{-}[dr] \ar@{.}@/^2.75pc/^{p_a^1}[rrr]&\ar@{}[r]|{w_1} & & \bullet\ar@{-}[ul]\ar@{--}[dl] \\
& \bullet\ar@{--}[r]\ar@{.}@/^.5pc/[ul]^{q_a^1} & a \\
& \bullet \ar@{--}[r]\ar@{-}[u] & a^* \\
\lck{\bullet}\ar@{-}[ur]\ar@{.}@/^.5pc/[ur]^{p_{a^*}^2}\ar@{.}@/_2.75pc/[rrr]_{q_{a^*}^2}\ar@{-}[dr] &\ar@{}[r]|{w_2} & & \bullet\ar@{--}[ul]\ar@{-}[dl]\ar@{-}[uuu] \\
& \bullet\ar@{-}[r] & \bullet \\
}}
+ (-1)^{d_a^2}
\raisebox{5pc}{
\xymatrix@R=1pc@C=1pc{
& \bullet \ar@{-}[r] & \bullet \\
*+[o][F-]{\bullet}\ar@{-}[ur]\ar@{-}[dr] \ar@{.}@/^2.75pc/^{p_{a}^2}[rrr]&\ar@{}[r]|{w_2} & & \bullet\ar@{-}[ul]\ar@{--}[dl] \\
& \bullet\ar@{--}[r]\ar@{.}@/^.5pc/[ul]^{q_{a}^2} & a \\
& \bullet \ar@{--}[r]\ar@{-}[u] & a^* \\
\lck{\bullet}\ar@{-}[ur]\ar@{.}@/^.5pc/[ur]^{p_{a^*}^1}\ar@{.}@/_2.75pc/[rrr]_{q_{a^*}^1}\ar@{-}[dr] &\ar@{}[r]|{w_1} & & \bullet\ar@{--}[ul]\ar@{-}[dl]\ar@{-}[uuu] \\
& \bullet\ar@{-}[r] & \bullet \\
}}
$$
\caption{the graded Kontsevich bracket $\{w_1,w_2\}$ in $\Neck_Q$. The dashed links, the beads inbetween these links (denoted as $a$ and $a^* = \frac{\partial}{\partial a}$ and the round ornamentation are removed and the open necklaces are connected as indicated. The exponents used are $d_a^1 = (|w_1|-1)|p_{a^*}^2|+|p_a^1||q_a^1|+1$ and $d_a^2 = (|w_2|-1)|p_{a}^2|+|p_{a^*}^1||q_{a^*}^1|$.}
\label{bracket}
\end{figure}

For a quiver $Q$, the necklace Lie algebra was introduced in \cite{RafLieven} in order to generalize the classical Karoubi-De Rham complex to noncommutative geometry. We will briefly recall the notions from \cite{RafLieven} needed for the remainder of this section.
\begin{defn}
For a quiver $Q$, define its \emph{double quiver} $\overline{Q}$ as the quiver obtained by adding for each arrow $a$ in $Q$ an arrow $a^*$ in the opposite direction of $a$ to $Q$.
\end{defn}
Now recall that the necklace Lie algebra was defined as 
\begin{defn}
The \emph{necklace Lie algebra} $N_Q$ is defined as
$N_Q := \CC\overline{Q}/[\CC\overline{Q},\CC\overline{Q}]$
equipped with the Kontsevich bracket which is defined on two necklaces $w_1$ and $w_2$ as illustrated in Figure \ref{bracket}. That is, for each arrow $a$ in $w_1$, look for all occurrences of $a^*$ in $w_2$, remove $a$ from $w_1$ and $a^*$ from $w_2$ and connect the corresponding open ends of both necklaces. Next, sum all the necklaces thus obtained. Now repeat the process with the roles of $w_1$ and $w_2$ reversed and deduct this sum from the first.
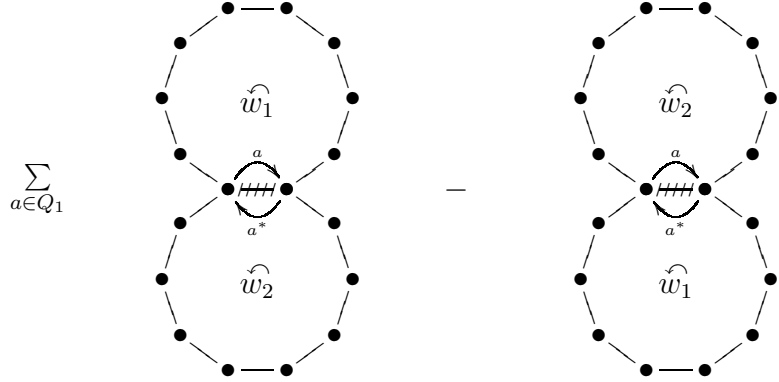
\begin{figure}
\[
\begin{xy}/r3pc/:
{\xypolygon10{~*{\bullet}~>>{}}},
"1" *+{\txt{\tiny $u$}},"8"="a","10"="a1","9"="a2","0"="c1",
"c1"+(0,-1.9),
{\xypolygon10{~*{\bullet}~>>{}}},
"5" *+{\txt{\tiny $v$}},"4"="b","6"="b1","5"="b2","0"="c2",
"a";"a2" **@{/};
"a" +(-2,0) *+{\txt{$\underset{a \in Q_1}{\sum}$}};
"c1" *+{\txt{$\overset{\curvearrowleft}{w_1}$}};
"c2" *+{\txt{$\overset{\curvearrowleft}{w_2}$}};
\POS"a" \ar@/^2ex/^{\txt{\tiny{$a$}}} "a2"
\POS"a2" \ar@/^2ex/^{\txt{\tiny{$a^*$}}} "a"
\end{xy}~\qquad
   \begin{xy}/r3pc/:
{\xypolygon10{~*{\bullet}~>>{}}},
"1" *+{\txt{\tiny $u$}},"8"="a","10"="a1","9"="a2","0"="c1",
"c1"+(0,-1.9),
{\xypolygon10{~*{\bullet}~>>{}}},
"5" *+{\txt{\tiny $v$}},"4"="b","6"="b1","5"="b2","0"="c2",
"a";"a2" **@{/};
"a"+(-2,0) *+{-};
"c1" *+{\txt{$\overset{\curvearrowleft}{w_2}$}};
"c2" *+{\txt{$\overset{\curvearrowleft}{w_1}$}};
\POS"a" \ar@/^2ex/^{\txt{\tiny{$a$}}} "a2"
\POS"a2" \ar@/^2ex/^{\txt{\tiny{$a^*$}}} "a"
\end{xy}
\]
\caption{Lie bracket $[ w_1,w_2 ]$ in $N_Q$.}
\label{bracket}
\end{figure}
\end{defn}
Now consider the following grading on $\CC\overline{Q}$: arrows $a$ in the original quiver are given degree $0$ and the starred arrows $a^*$ in $\overline{Q}$ are given degree $1$.  We now can consider the \emph{graded necklace Lie algebra} $\CC\overline{Q}/[\CC\overline{Q},\CC\overline{Q}]_{super}$ equipped with a graded version of the Kontsevich bracket, as introduced in \cite{lazaroiu}.
\begin{defn}
The graded necklace Lie algebra is defined as
$$\Neck_Q := \CC\overline{Q}/[\CC\overline{Q},\CC\overline{Q}]_{super}$$
equipped with the \emph{graded Kontsevich bracket} defined in Figure \ref{gradedKbracket}.
Monomials in $\Neck_Q$ are depicted as \emph{ornate necklaces}\index{ornate necklace}, where the beads represent arrows in the necklace and where one bead is encased, indicating the starting point of the necklace. 
\end{defn}
An example of an ornate necklace is
$$
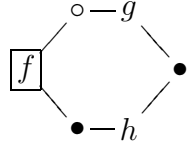
\xymatrix@R=.75pc@C=.75pc{
& \circ \ar@{-}[r] & g \\
\lck{f}\ar@{-}[ur]\ar@{-}[dr] & & & \bullet\ar@{-}[ul]\ar@{-}[dl] \\
& \bullet\ar@{-}[r] & h
}
$$
representing the element $f\delta g \Delta h \Delta$ if we let $\circ$ represent $\delta$ and $\bullet$ represents $\Delta$. The identities coming from dividing out supercommutators then look like
$$
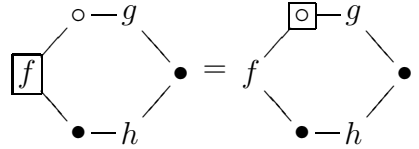
\raisebox{1.95pc}{\xymatrix@R=.75pc@C=.75pc{
& \circ \ar@{-}[r] & g \\
\lck{f}\ar@{-}[ur]\ar@{-}[dr] & & & \bullet\ar@{-}[ul]\ar@{-}[dl] \\
& \bullet\ar@{-}[r] & h
}}
= 
\raisebox{1.95pc}{\xymatrix@R=.75pc@C=.75pc{
& \lck{\circ} \ar@{-}[r] & g \\
f\ar@{-}[ur]\ar@{-}[dr] & & & \bullet\ar@{-}[ul]\ar@{-}[dl] \\
& \bullet\ar@{-}[r] & h
}}.
$$

This graded necklace Lie algebra is the noncommutative equivalent of the classical graded Lie algebra of multivector fields.
\begin{theorem}\label{MultiIsNecklace}
Let $Q$ be a quiver, then
$$T_{\CC Q}\DDer(\CC Q)/[T_{\CC Q}\DDer(\CC Q),T_{\CC Q}\DDer(\CC Q)] \cong \Neck_Q$$
as graded Lie algebras.
\end{theorem}
\begin{proof}
From \cite{MichelDPA} we know that the module of double derivations $\DDer(\CC Q)$ is generated as a $\CC Q$-bimodule by the double derivations $\frac{\partial}{\partial a}$, $a\in Q_1$, defined as
$$\frac{\partial}{\partial a}(b) = \left\{\begin{array}{ll} e_{t(a)}\otimes e_{h(a)} &\, b = a \\ 0 &\, b\neq a\end{array}\right..$$
Now note that we may identify $\frac{\partial}{\partial a}$ with an arrow $a^*$ in $\overline{Q}$ in the opposite direction of $a$: $\frac{\partial}{\partial a} = e_{h(a)}\frac{\partial}{\partial a} e_{t(a)}$, so $T_{\CC Q}\DDer(\CC Q) \cong \CC\overline{Q}$. The arrows $a^*$ correspond to the degree $1$ elements in the tensor algebra and the original arrows to the degree $0$ arrows. That is, supercommutators in the algebra on the left correspond to supercommutators in the algebra on the right.

Now note that the NC Schouten bracket on a path in $\overline{Q}$ becomes
\begin{eqnarray*}
\db{a^*,a_1.\dots.a_n} & = & \sum_{i=0}^{n-1} (-1)^{|a_1|+\dots+|a_i|}a_1\dots a_i\db{a^*,a_{i+1}}a_{i+2}\dots a_n \\
& = & \sum_{i=0, a_{i+1} = a}^{n} (-1)^{|a_1|+\dots+|a_i|}a_1\dots a_i\otimes a_{i+2}\dots a_n,
\end{eqnarray*}
where $a_0 = {t(a_1)}$.
But this is the graded version of the necklace Loday algebra considered in \cite{MichelDPA}. This becomes the graded necklace Lie algebra when restricting to closed paths and modding out commutators.
\end{proof}

An immediate corollary of the Theorem above is
\begin{cor}\label{DPT}
A nonzero double bracket on the path algebra $\CC Q$ is completely determined by a linear combination of necklaces of degree $2$.
\end{cor}
For the remainder of the paper, we will assume the brackets to be nonzero.

\section{Linear Double Poisson Structures}\label{DoubleLiePoisson}
In classical Poisson geometry, linear Poisson structures are defined on $\CC^n$ through a Poisson tensor of the form
$$
\pi = c_{ij}^k\,x_k\,\frac{\partial}{\partial
  x_i}\wedge\frac{\partial}{\partial x_j},
$$
where we use Einstein notation, that is, we sum over repeated indices. For this expression to be a Poisson tensor, the constant factors $c_{ij}^k$ must satisfy
$$c_{jk}^hc_{hi}^p+c_{ki}^hc_{hj}^p+c_{ij}^hc_{hk}^p = 0,$$
i.e. $c_{ij}^k$ are the structure constants of an $n$-dimensional Lie algebra.

In order to translate this setting to NC Poisson geometry, we first of all note that the role of affine space is assumed by the representation spaces of quivers and the bivector $\frac{\partial}{\partial x_i}\wedge\frac{\partial}{\partial x_j}$ is replaced by the degree $2$ part of a necklace, $\frac{\partial}{\partial a}\frac{\partial}{\partial b}$, by Corollary \ref{DPT}. We define
\begin{defn}
Let $Q$ be a quiver.
A \emph{linear double bracket} on $\CC Q$ is a double bracket determined by a double tensor of the form
$$P_{lin} := \sum_{\stackrel{a,b,c\in Q_1}{a\frac{\partial}{\partial b}\frac{\partial}{\partial c}\neq 0\in\mathfrak{n}_Q}} c^a_{bc}a\frac{\partial}{\partial b}\frac{\partial}{\partial c},$$
with all $c_{bc}^a \in \CC$.
\end{defn}
We can characterize the linear double Poisson brackets as follows.
\begin{theorem}\label{linbrack}
A linear double bracket 
$$P_{lin} := \sum_{\stackrel{a,b,c\in Q_1}{a\frac{\partial}{\partial b}\frac{\partial}{\partial c}\neq 0\in\mathfrak{n}_Q}} c^a_{bc}a\frac{\partial}{\partial b}\frac{\partial}{\partial c},$$
Poisson bracket if and only if for all $p,q,r,s\in Q_1$ such that $p\frac{\partial}{\partial q}\frac{\partial}{\partial r}\frac{\partial}{\partial s}\neq 0 \in \mathfrak{n}_Q$ we have
$$\sum_{x\in Q_1^{(p,q,rs)}}c_{rs}^xc_{qx}^p - \sum_{y\in Q_1^{(qr,p,s)}}c_{ys}^pc_{qr}^y,$$
where 
$$Q_1^{(p,q,rs)} = \{a\in Q_1 \mid a\frac{\partial}{\partial r}\frac{\partial}{\partial s}, p\frac{\partial}{\partial q}\frac{\partial}{\partial a} \neq 0 \in \mathfrak{n}_Q\}$$
and
$$Q_1^{(qr,p,s)} = \{a\in Q_1 \mid a\frac{\partial}{\partial q}\frac{\partial}{\partial r}, p\frac{\partial}{\partial a}\frac{\partial}{\partial s} \neq 0 \in \mathfrak{n}_Q\}.$$
\end{theorem}
\begin{proof}
We have to verify when $\{P_{lin},P_{lin}\} = 0$ modulo commutators. First of all observe that a straightforward computation yields that for any $x,y,z,u,v,w\in Q_1$
\begin{eqnarray*}
\{x\frac{\partial}{\partial y}\frac{\partial}{\partial z},u\frac{\partial}{\partial v}\frac{\partial}{\partial w}\} & = &
\delta_{xw} u \frac{\partial}{\partial v}\frac{\partial}{\partial y}\frac{\partial}{\partial z}
+ \delta_{uz} x \frac{\partial}{\partial y}\frac{\partial}{\partial v}\frac{\partial}{\partial w}\\
& & - \delta_{uy} x \frac{\partial}{\partial v}\frac{\partial}{\partial w}\frac{\partial}{\partial z}
- \delta_{xv} u \frac{\partial}{\partial y}\frac{\partial}{\partial z}\frac{\partial}{\partial w}
\end{eqnarray*}
modulo commutators. Next, note that this equality implies that $\{P_{lin},P_{lin}\}$ lies in the sub vector space of $\mathfrak{n}_Q$ that has as basis $B$ all ornate necklaces of the form $p\frac{\partial}{\partial q}\frac{\partial}{\partial r}\frac{\partial}{\partial s}$ with $p,q,r,s\in Q_1$ and where $p$ is the encased bead. We now write
\begin{eqnarray*}
\{P_{lin},P_{lin}\} & = & \sum_{x,y,z}\sum_{u,v,w} \{x\frac{\partial}{\partial y}\frac{\partial}{\partial z},u\frac{\partial}{\partial v}\frac{\partial}{\partial w}\} \\
&=& \sum_{x,y,z}\left ( \sum_{u,v} c_{yz}^xc_{vx}^u u \frac{\partial}{\partial v}\frac{\partial}{\partial y}\frac{\partial}{\partial z}
+  \sum_{v,w} c_{yz}^xc_{vw}^z x \frac{\partial}{\partial y}\frac{\partial}{\partial v}\frac{\partial}{\partial w}\right. \\
& & \left. -  \sum_{v,w} c_{yz}^xc_{vw}^y x \frac{\partial}{\partial v}\frac{\partial}{\partial w}\frac{\partial}{\partial z} -  \sum_{u,w} c_{yz}^xc_{xw}^u u \frac{\partial}{\partial y}\frac{\partial}{\partial z}\frac{\partial}{\partial w} \right),
\end{eqnarray*}
where in order to lighten notation, we do not explicitly write down the additional restrictions on $x,y,z,u,v,w$ for a $c_{yz}^x$ and $c_{uv}^w$ to be defined. Regrouping this expression with respect to the basis $p\frac{\partial}{\partial q}\frac{\partial}{\partial r}\frac{\partial}{\partial s}$, we get
$$
\{P_{lin},P_{lin}\} = \sum_{p\frac{\partial}{\partial q}\frac{\partial}{\partial r}\frac{\partial}{\partial s}\in B} 2\left(\sum_{x\in Q_1^{(p,q,rs)}}c_{rs}^xc_{qx}^p - \sum_{y\in Q_1^{(qr,p,s)}}c_{ys}^pc_{qr}^y\right)p\frac{\partial}{\partial q}\frac{\partial}{\partial r}\frac{\partial}{\partial s},
$$
where the summation sets simply list which coefficients are defined. Equating this last expression to zero concludes the proof.
\end{proof}

We have the following immediate corollary for the induced Poisson bracket on the representation and quotient spaces.
\begin{cor}
Let $P_{lin}$ be a linear double Poisson bracket, then the induced bracket $tr(P_{lin})$ on $\rep_n(\CC Q)$ and on $\iss_n(\CC Q)$ is a linear Poisson bracket.
\end{cor}
\begin{proof}
This is immediate from the definition of the induced bracket.
\end{proof}

Note that this proof immediately indicates that not all linear Poisson structures are induced by linear double Poisson structures, as the condition on the coefficients of the latter is more restrictive than the condition on the coefficients of the former. We even have
\begin{cor}\label{assalg}
With notations as above, let $P_{lin}$ be a linear double Poisson bracket and $v\in Q_0$, then $P_{lin}$ induces an associative algebra structure on the vector space generated by all loops in $v$.
\end{cor}
\begin{proof}
It suffices to observe that the condition in Theorem \ref{linbrack} exactly determines the structure constants of an associative algebra structure on the loops in a vertex.
\end{proof}
\section{Cohomology of linear double Poisson brackets and Hochschild cohomology}\label{example}
In this section, we will first of all specialize the description of
linear double Poisson brackets obtained in the previous section to the
free algebra in $n$ variables. Then, we will give a link between the
double Poisson cohomology of such linear double Poisson brackets and
the Hochschild cohomology of the corresponding associative algebra.

First of all observe that Corollary \ref{assalg} becomes much stronger for the free algebra in $n$ variables.
\begin{prop}\label{dpass1to1}
There is a one-to-one correspondence between linear double Poisson brackets on $\CC\langle x_1,\dots,x_n\rangle$ and associative algebra structures on $V = \CC x_1\oplus\dots\oplus\CC x_n$. Explicitly, consider the associative algebra structure on $V$ determined by 
$$x_ix_j := \sum_{i,j,k=1}^n c_{ij}^k x_k,$$
then the corresponding double Poisson bracket is given by
$$\db{x_i,x_j} = \sum_{i,j,k=1}^n (c_{ij}^k x_k\otimes 1 - c_{ji}^k 1 \otimes x_k).$$
\end{prop}
\begin{proof}
This holds as the first condition in Theorem \ref{linbrack} vanishes because the free algebra in $n$ variables is the path algebra of a quiver with one vertex and $n$ loops.
\end{proof}
As a corollary we obtain for the noncommutative affine plane
\begin{cor}
Up to affine transformation, there are only $7$ different linear double Poisson brackets on $\Cxy$. Their double Poisson brackets are
$$
\begin{array}{lcl}
P_{lin}^{\CC\times\CC} &  = & x\frac{\partial}{\partial x}\frac{\partial}{\partial x} + y\frac{\partial}{\partial y}\frac{\partial}{\partial y} \\
P_{lin}^{\CC\times\CC\epsilon^2} &  = & x\frac{\partial}{\partial x}\frac{\partial}{\partial x} \\
P_{lin}^{\CC\oplus\CC\epsilon^2} &  = & x\frac{\partial}{\partial x}\frac{\partial}{\partial x} + y\frac{\partial}{\partial x}\frac{\partial}{\partial y} + y \frac{\partial}{\partial y}\frac{\partial}{\partial x} \\
P_{lin}^{\CC\epsilon\oplus\CC\epsilon^2} & = & y\frac{\partial}{\partial x}\frac{\partial}{\partial x} \\
P_{lin}^{B_2^1} & = & x\frac{\partial}{\partial x}\frac{\partial}{\partial x} + y\frac{\partial}{\partial x}\frac{\partial}{\partial y}\\
P_{lin}^{B_2^2} & = & x\frac{\partial}{\partial x}\frac{\partial}{\partial y} + y\frac{\partial}{\partial y}\frac{\partial}{\partial y}\\
P_{lin}^{\CC\epsilon^2\oplus\CC\epsilon^2} & = & 0
\end{array}
$$
\end{cor}
\begin{proof}
This follows immediately from Proposition \ref{dpass1to1} and the classification of all (non-unital) $2$-dimensional associative algebras obtained in \cite{gab}. The upper indices of the brackets listed here correspond to the algebra structures they induce.
\end{proof}
Moreover, there is a direct connection between the Hochschild cohomology of the finite dimensional algebra and the double Poisson cohomology of the free algebra.
\begin{theorem}\label{hochschild}
Let $A = \CC x_1 \oplus \dots \oplus \CC x_n$ be an $n$-dimensional vector space and let 
$$P = \sum_{i,j,k=1}^n c_{ij}^kx_k\frac{\partial}{\partial x_i}\frac{\partial}{\partial x_j}$$
be a linear double Poisson structure on $T_\CC A = \CC\langle x_1,\dots, x_n\rangle$. Consider $A$ as an algebra through the product induced by the structure constants of $P$ and let $HH^\bullet(A)$ denote the Hochschild cohomology of this algebra, then
$$(H_P^\bullet(T_\CC A))_1 \cong HH^\bullet(A).$$
Here the grading on $(H_P^\bullet(T_\CC A))$ is induced by the grading on $\left (\frac{T_{T_\CC A}}{[T_{T_\CC A},T_{T_\CC A}]}\right )_{i}$, which is defined through $\deg(x_i)  =1$.
\end{theorem}
\begin{proof}
First of all observe that we have a basis for $\left (\frac{T_{T_\CC A}}{[T_{T_\CC A},T_{T_\CC A}]}\right )_{i,1}$ consisting of all possible elements of the form
$$\frac{\partial}{\partial x_{k_1}}\frac{\partial}{\partial x_{k_2}}\dots\frac{\partial}{\partial x_{k_i}}x_{\ell}.$$
Now consider the linear map
$$\varphi_i: \underbrace{A^*\otimes A^*\otimes\dots \otimes A^*}_{i \mathrm{~factors}} \otimes A\rightarrow
\left (\frac{T_{T_\CC A}}{[T_{T_\CC A},T_{T_\CC A}]}\right )_{i,1}$$
defined as 
$$\varphi_{i}(x^*_{k_1}\otimes x^*_{k_2}\otimes \dots \otimes x^*_{k_i}\otimes x_\ell) = \frac{\partial}{\partial x_{k_1}}\frac{\partial}{\partial x_{k_2}}\dots\frac{\partial}{\partial x_{k_i}}x_{\ell},$$
then $\varphi := (\varphi_i)$ is a morphism of complexes
$$\xymatrix{
0\ar[r] & A\ar[r]^d\ar[d]^{\varphi_0} & A^*\otimes A\ar[r]^d\ar[d]^{\varphi_1} & A^*\otimes A^*\otimes A\ar[r]\ar[d]^{\varphi_2} & \dots \\
0\ar[r] & \left (\frac{T_{T_\CC A}}{[T_{T_\CC A},T_{T_\CC A}]}\right )_{0,1}\ar[r]^{d_P} &\left  (\frac{T_{T_\CC A}}{[T_{T_\CC A},T_{T_\CC A}]}\right )_{1,1}\ar[r]^{d_P} &\left  (\frac{T_{T_\CC A}}{[T_{T_\CC A},T_{T_\CC A}]}\right )_{2,1}\ar[r] & \dots
}$$
where the upper complex is the Hochschild complex with $d$ defined as
\begin{align*}
d(x^*_{k_1}\otimes x^*_{k_2}\otimes \dots \otimes x^*_{k_i}\otimes x_\ell) = &
\sum_{s,t=1}^n c_{s\ell}^t\, x_{s}^*\otimes x^*_{k_1}\otimes x^*_{k_2}\otimes \dots \otimes x^*_{k_i}\otimes x_t \\
+&\sum_{r = 1}^{i-1}(-1)^r\sum_{s,t=1}^n c_{st}^{k_r}\,x^*_{k_1}\otimes \dots \otimes \underbrace{x^*_s\otimes x^*_t}_{r\mathrm{-th~factor}}\otimes\dots \otimes x^*_{k_i}\otimes x_\ell \\
+& (-1)^{i+1}\sum_{s,t=1}^n c_{\ell s}^t \,x^*_{k_1}\otimes x^*_{k_2}\otimes \dots \otimes x^*_{k_i}\otimes x_{s}^*\otimes  x_t.
\end{align*}
In order to see $\varphi$ is a morphism of complexes, we compute
\begin{eqnarray*}
\db{P,\frac{\partial}{\partial x_{k_1}}\frac{\partial}{\partial
    x_{k_2}}\dots\frac{\partial}{\partial x_{k_i}}x_{\ell}} & = &
\underbrace{(-1)^i\frac{\partial}{\partial
    x_{k_1}}\dots\frac{\partial}{\partial x_{k_i}}\db{P,x_\ell}}_{T_1}
\\
&&+ \underbrace{\db{P,\frac{\partial}{\partial x_{k_1}}\frac{\partial}{\partial x_{k_2}}\dots\frac{\partial}{\partial x_{k_i}}}x_\ell}_{T_2}.
\end{eqnarray*}
The first term in this expression is mapped by the multiplication on the tensor algebra to
\begin{eqnarray*}
\mu(T_1) &=& (-1)^i
\frac{\partial}{\partial x_{k_1}}\frac{\partial}{\partial x_{k_2}}\dots\frac{\partial}{\partial x_{k_i}}\{P,x_\ell\}\\
& = & (-1)^i\frac{\partial}{\partial x_{k_1}}\frac{\partial}{\partial x_{k_2}}\dots\frac{\partial}{\partial x_{k_i}}\sum_{r,s=1}^n\left ( c_{s\ell}^rx_r\frac{\partial}{\partial x_s} - c_{\ell s}^r\frac{\partial}{\partial x_s}x_r\right ) \\
& = & \sum_{r,s=1}^n c_{s\ell}^r \frac{\partial}{\partial x_s}\frac{\partial}{\partial x_{k_1}}\frac{\partial}{\partial x_{k_2}}\dots\frac{\partial}{\partial x_{k_i}}x_r \\
& & -(-1)^i\sum_{r,s=1}^n c_{\ell s}^r\frac{\partial}{\partial x_{k_1}}\frac{\partial}{\partial x_{k_2}}\dots\frac{\partial}{\partial x_{k_i}}\frac{\partial}{\partial x_s}x_r
\end{eqnarray*}
The second term in the expression is mapped to
\begin{eqnarray*}
\mu(T_2) & = & \sum_{r=0}^{i-1} (-1)^r \frac{\partial}{\partial x_{k_1}}\dots\frac{\partial}{\partial x_{k_r}}\{P,\frac{\partial}{\partial x_{k_{r+1}}}\}\frac{\partial}{\partial x_{k_{r+2}}}\dots\frac{\partial}{\partial x_{k_i}}x_\ell \\
& = & -\sum_{r=0}^{i-1} (-1)^r \sum_{s,t=1}^n c_{st}^{k_{r+1}}\frac{\partial}{\partial x_{k_1}}\dots\frac{\partial}{\partial x_{k_r}}
\frac{\partial}{\partial x_s}\frac{\partial}{\partial x_t}
\frac{\partial}{\partial x_{k_{r+2}}}\dots\frac{\partial}{\partial x_{k_i}}x_\ell.
\end{eqnarray*}
Adding these two expressions together indeed yields
\begin{eqnarray*}
d_P\left(\frac{\partial}{\partial x_{k_1}}\frac{\partial}{\partial
  x_{k_2}}\dots\frac{\partial}{\partial x_{k_i}}x_{\ell}\right) &=&
d_P(\varphi_i(x^*_{k_1}\otimes x^*_{k_2}\otimes \dots \otimes x^*_{k_i}\otimes x_\ell))\\
&=& \varphi_i(d(x^*_{k_1}\otimes x^*_{k_2}\otimes \dots \otimes x^*_{k_i}\otimes x_\ell)),
\end{eqnarray*}
so we have a morphism of complexes. It is easy to see this is an isomorphism when restricting to degree $1$ terms in the lower complex, which finishes the proof.
\end{proof}

\begin{remark}
Note that Theorem \ref{hochschild} can be seen as the noncommutative analogue of the fact that for a Lie-Poisson structure associated to a compact Lie group the Poisson cohomology can be written as the tensor product of the Lie algebra cohomology with the Casimir elements of the Poisson bracket.
\end{remark}

The relation between higher degree components of the double Poisson cohomology and the Hochschild cohomology of the corresponding finite dimensional algebra is less obvious. However, we have
\begin{theorem}\label{Casimirs}
With notations as in the previous theorem, there is a canonical embedding
$$H^0_P(T_\CC A)\hookrightarrow \CC\oplus\bigoplus_{k\geq 1} HH^0(A,A^{\otimes k}),$$
where the action of $A$ on $A^{\otimes k}$ is the inner action on the two outmost copies of $A$ in the tensor product.
\end{theorem}
\begin{proof}
We have already shown in Theorem \ref{hochschild} that the Hochschild cohomology $HH^0(A)$ corresponds to the degree $1$ terms in $H_P^0(T_\CC A)$. We will now show that
$(H_P^0(T_\CC A))_k \hookrightarrow HH^0(A,A^{\otimes k})$. Consider the map
$$\varphi_0: \frac{T_\CC A}{[T_\CC A,T_\CC A]}\rightarrow A^{\otimes k}:x_{i_1}\otimes\dots\otimes x_{i_k} \mapsto \sum_{\ell=0}^{k-1} \sigma_{(1\dots k)}^\ell x_{i_1}\otimes\dots\otimes x_{i_k},$$
where for $s\in S_k$ a permutation we define $\sigma_s(a_1\otimes\dots\otimes a_k) := a_{s(1)}\otimes\dots\otimes a_{s(k)}$. It is easy to see that this map is well-defined. We will also need the map
$$\varphi_1: \left (\frac{T_{T_\CC A}}{[T_{T_\CC A},T_{T_\CC A}]}\right )_{1,k}\rightarrow A^*\otimes A^{\otimes k}$$
defined as
$$\varphi_1(\frac{\partial}{\partial x_i}x_{j_1}\dots x_{j_k}) = x_i^*\otimes x_{j_1}\otimes\dots\otimes x_{j_k},$$
where we fixed a similar basis as in Theorem \ref{hochschild} for $\left (\frac{T_{T_\CC A}}{[T_{T_\CC A},T_{T_\CC A}]}\right )_{1,k}$.

Now compute for $f = x_{i_1}\dots x_{i_k}$ (using Einstein notations)
\begin{eqnarray*}
d_P(f) &=& 
\left(c_{pq}^r\,\frac{\partial f}{\partial x_q}''\frac{\partial f}{\partial x_q}'x_r - c_{qp}^rx_r\frac{\partial f}{\partial x_q}''\frac{\partial f}{\partial x_q}'\right)\frac{\partial}{\partial x_p} \\
& = &\sum_{s=0}^{k-1} \left(c_{pq}^r\,\frac{\partial x_{i_{s+1}}}{\partial x_q}''x_{i_{s+2}}\dots x_{i_k}x_{i_1}\dots x_{i_s}\frac{\partial x_{i_{s+1}}}{\partial x_q}'x_r\right.\\
 &  & \left. - c_{qp}^r\,x_r\frac{\partial x_{i_{s+1}}}{\partial x_q}''x_{i_{s+2}}\dots x_{i_k}x_{i_1}\dots x_{i_s}\frac{\partial x_{i_{s+1}}}{\partial x_q}'\right)\frac{\partial}{\partial x_p} \\
& = &  \sum_{s=0}^{k-1} \left(c_{pi_{s+1}}^r\,x_{i_{s+2}}\dots x_{i_k}x_{i_1}\dots x_{i_s}x_r\right.\\
 &  & \left. - c_{i_{s+1}p}^rx_rx_{i_{s+2}}\dots x_{i_k}x_{i_1}\dots x_{i_s}\right)\frac{\partial}{\partial x_p}
\end{eqnarray*}
This is mapped by $\varphi_1$ to
\begin{align*}
 &\sum_{s=0}^{k-1} (c_{pi_{s+1}}^r\,x_p^*\otimes x_{i_{s+2}}\otimes\dots\otimes x_{i_k}\otimes x_{i_1}\otimes\dots\otimes x_{i_s}\otimes x_r \\
 -& c_{i_{s+1}p}^r\,x_p^*\otimes x_r\otimes x_{i_{s+2}}\otimes\dots\otimes x_{i_k}\otimes x_{i_1}\otimes\dots\otimes x_{i_s}).
\end{align*}
On the other hand, we have 
\begin{eqnarray*}
d\varphi_0(f) & = & \sum_{s=0}^{k-1}  d(x_{i_{s+1}}\otimes\dots\otimes x_{i_k}\otimes x_{i_1}\otimes\dots\otimes x_{i_s}) \\
& = & \sum_{s=0}^{k-1}\left( c^r_{pi_{s}}\,x_p^* \otimes  x_{i_{s+1}}\otimes\dots\otimes x_{i_k}\otimes x_{i_1}\otimes\dots\otimes x_{i_{s-1}}\otimes x_r\right. \\
& &\left. - c_{i_{s+1}p}^r\,x_p^*\otimes x_r\otimes x_{i_{s+2}}\otimes\dots\otimes x_{i_k}\otimes x_{i_1}\otimes\dots\otimes x_{i_s}\right).
\end{eqnarray*}
So, after reindexing the first term, we obtain $d(\varphi_0(f)) = \varphi_1(d_P(f))$, which finishes the proof.
\end{proof}

\section{Examples of $H^0$ and $H^1$ for several double Poisson structures on $\Cxy$}\label{H0H1}

In this section, we will determine the double Poisson cohomology
groups $H^0$ and $H^1$, of three different double Poisson structures
on the free algebra $\Cxy$. To do this, we will use some tools and
techniques that could be useful to compute other double Poisson
cohomology groups, associated to other double Poisson algebras.

In order to compute the double Poisson cohomology of $\Cxy$ with the
different brackets we will introduce below, the following noncommutative version of the Euler formula will prove useful.
\begin{prop}[NC-Euler formula]\label{nceuler}
Let $Q$ be a quiver, $p$ a path in $Q$ and $a$ an arrow of $Q$, then
$$\mu\circ\left(a\frac{\partial}{\partial a}\right)(p) = (\deg_a(p))\,p.$$
Whence
$$\sum_{a\in Q_1} \mu\circ \left(a\frac{\partial}{\partial a}\right)(p) = |p|\,p,$$
with $|p|$ the length of the necklace.
\end{prop}
\begin{proof}
From the definition of $a\frac{d}{da}$, for any path in $Q$, of the form
$a_1\cdot a_2\cdots a_{n-1}\cdot a_n$ (with $n\in\N^*$ and
$a_1,a_2,\dots,a_n$, some arrows of $Q$), we have:
\begin{eqnarray*}
\frac{\partial}{\partial a}\left(a_1\cdot a_2\cdots a_{n-1}\cdot a_n\right) 
&=& \sum_{i=1}^n  a_1\cdots a_{i-1}\cdot\frac{\partial a_i}{\partial
  a}\cdot a_{i+1}\cdots a_n\\
&=& \sum_{\stackrel{i=1}{a=a_i}}^n  
a_1\cdots a_{i-1}\cdot e_{t(a_i)}\otimes e_{h(a_i)}\cdot a_{i+1}\cdots a_n\\
&=& \sum_{\stackrel{i=1}{a=a_i}}^n  
a_1\cdots a_{i-1}\otimes a_{i+1}\cdots a_n.
\end{eqnarray*}
So that, by definition of the inner product, we obtain:
\begin{eqnarray*}
\left(a\frac{\partial}{\partial a}\right)\left(a_1\cdot a_2\cdots a_{n-1}\cdot a_n\right) 
&=& \sum_{\stackrel{i=1}{a=a_i}}^n  
a_1\cdots a_{i-1}\otimes a_i\cdot a_{i+1}\cdots a_n,
\end{eqnarray*}
proving the proposition.
\end{proof}

The rest of this section will be devoted to the determination of the
low-dimensional double Poisson cohomology groups $H^0$ and $H^1$ of
the free algebra $\Cxy$, equipped with the following (linear and non-linear) double
Poisson tensors:
\begin{enumerate}
\item the linear double Poisson brackets $P_0 = x\frac{d}{dx}\frac{d}{dx}$ and 
$x\frac{d}{dx}\frac{d}{dx}+y\frac{d}{dy}\frac{d}{dy}$, for which the
 corresponding Poisson bracket on $\rep_1(\Cxy)=\CC[x,y]$ is zero;
\item the linear double Poisson brackets 
$P_1=x \frac{d}{dx}\frac{d}{dx} + y\frac{d}{dx}\frac{d}{dy}$ or 
$x \frac{d}{dx}\frac{d}{dy} + y\frac{d}{dy}\frac{d}{dy}$, for which
  the Poisson bracket obtained with the trace on $\rep_1(\Cxy)$ is a
  non-trivial Poisson bracket;
\item the quadratic double Poisson bracket
  $P=x\frac{d}{dx}x\frac{d}{dy}$, which induces a quadratic
  (non-trivial) Poisson bracket on $\rep_1(\Cxy)$.
\end{enumerate}
\subsection{The linear double Poisson tensors
  $x\frac{d}{dx}\frac{d}{dx}$ and $x\frac{d}{dx}\frac{d}{dx}+y\frac{d}{dy}\frac{d}{dy}$}

Let us first consider the linear double Poisson tensor, given by 
$$
P_0:=P_{lin}^{\CC\times\CC\epsilon^2}= x \frac{d}{dx}\frac{d}{dx}.
$$
Our aim, in this subsection, is to give an explicit basis of the double Poisson
cohomology groups $H^0_{P_0}(\Cxy)$ and
$H^1_{P_0}(\Cxy)$, associated to this double Poisson tensor $P_0$.

First of all, let us consider the operator $d^0_{P_0}$ and the space $H^0_{P_0}(\Cxy)$. 
\begin{prop}\label{imageofdp00}
For $f\in\Cxy$, we have
$$d_{P_0}(f) = 
\xymatrix@R=1.75pc{\lck{\circ}\ar@{-}@/^.25pc/[r] & 
\left(\frac{df}{dx}''\frac{df}{dx}'x -\frac{df}{dx}''\frac{df}{dx}'x\right) \ar@{-}@/^.5pc/[l]}.
$$
Which leads to
$$H_{P_0}^0(\Cxy) \simeq \CC[x]\oplus\CC[y].$$
\end{prop}
\begin{proof}
Let $f\in\Cxy$ and let recall that $P_0$ denotes the double Poisson
tensor $P_0=x \frac{d}{dx}\frac{d}{dx}$. Then, using the properties of
the double Gerstenhaber bracket $\db{-,-}$, given in \cite[\S
  2.7]{MichelDPA}, we compute the
double Schouten-Nijenhuis bracket of $f$ and $P_0$:
\begin{eqnarray*}
\db{f,P_0} &=& \db{f,x \frac{d}{dx}\frac{d}{dx}}\\
   &=& -x \frac{d}{dx}\db{f,\frac{d}{dx}}+
     \db{f,x\frac{d}{dx}}\frac{d}{dx}\\
   &=& x \frac{d}{dx}\frac{df}{dx}''\otimes \frac{df}{dx}'
     -x\frac{df}{dx}''\otimes \frac{df}{dx}'\frac{d}{dx},
\end{eqnarray*}
so that, computing modulo the commutators, we obtain exactly
$$
d_{P_0}(f) = \mu\left(\db{P_0,f}\right) = 
\left( \frac{df}{dx}''\frac{df}{dx}'x
     -x\frac{df}{dx}'' \frac{df}{dx}'\right)\frac{d}{dx}.
$$
Then, a $0$-cocycle for the double Poisson cohomology, corresponding
to $P_0$ is an element $f\in\Cxy$ satisfying $d_{P_0}(f) = 0$,
which means that 
$$
\frac{df}{dx}''\frac{df}{dx}'x
     -x\frac{df}{dx}'' \frac{df}{dx}' = 0,
$$
that is to say, the element $\displaystyle\frac{df}{dx}'' \frac{df}{dx}'\in\Cxy$
commutes with $x$, so is necessarily an element of $\CC[x]$. According to the
NC-Euler Formula (Proposition \ref{nceuler}), we have 
$$
\deg_x(f)\, f = \mu\circ\left(x\frac{\partial}{\partial x}\right)(f) 
= \frac{df}{dx}'x \frac{df}{dx}'' \in \frac{df}{dx}'' \frac{df}{dx}'x + \lbrack\Cxy,\Cxy\rbrack,
$$
so that, modulo commutators, we either have $\deg_x(f)=0$ and hence $f\in\CC[y]$,
or $f\in\CC[x]$. But then
$$
H_{P_0}^0(\Cxy) \simeq \CC[x]\oplus\CC[y].
$$
\end{proof}
Next, let us determine the first double Poisson cohomology group,
related to $P_0$. First of all, observe that an
element of $(T_{\Cxy}/[T_{\Cxy},T_{\Cxy}])_1$ can be uniquely written as 
$\displaystyle f\frac{d}{dx} + g\frac{d}{dy}$, with $f,g\in\Cxy$.
By a direct computation (analogous to what we did for $d_{P_0}(f)$), 
we obtain the value of the coboundary
operator $d_{P_0}$ on such an element. We obtain that
$$
d_{P_0}\left(f\frac{d}{dx} + g\frac{d}{dy}\right) = \Phi_1(f) + \Phi_2(g),
$$
where the operators $\Phi_1$ and $\Phi_2$ from $\Cxy$ to 
$(T_{\Cxy}/[T_{\Cxy},T_{\Cxy}])_2$ are defined,
for $f,g\in\Cxy$, by: 
\begin{align}\label{def:Phi_1}
\Phi_1(f) := & 
-\,\raisebox{1.8pc}{\xymatrix@R=.75pc@C=.75pc{& \circ & \\ 
\lck{f}\ar@{-}[ur]\ar@{-}[dr] & & 1\ar@{-}[ul]\ar@{-}[dl]\\ & \circ & }}
+\,\raisebox{1.8pc}{\xymatrix@R=.75pc@C=.75pc{& \circ & \\ 
\lck{\frac{df}{dx}'x}\ar@{-}[ur]\ar@{-}[dr] & & \frac{df}{dx}''
\ar@{-}[ul]\ar@{-}[dl]\\ & \circ & }}
-\,\raisebox{1.8pc}{\xymatrix@R=.75pc@C=.75pc{& \circ & \\ 
\lck{\frac{df}{dx}'}\ar@{-}[ur]\ar@{-}[dr] & & x\frac{df}{dx}''
\ar@{-}[ul]\ar@{-}[dl]\\ & \circ & }}\,,
\end{align}
and
\begin{align}\label{def:Phi_2}
\Phi_2(g) := &\, \raisebox{1.8pc}{\xymatrix@R=.75pc@C=.75pc{& \circ & \\ 
\lck{\frac{dg}{dx}'x}\ar@{-}[ur]\ar@{-}[dr] & & \frac{dg}{dx}''
\ar@{-}[ul]\ar@{-}[dl]\\ & \bullet & }}  
- \raisebox{1.8pc}{\xymatrix@R=.75pc@C=.75pc{& \circ & \\ 
\lck{\frac{dg}{dx}'}\ar@{-}[ur]\ar@{-}[dr] & & x\frac{dg}{dx}''
\ar@{-}[ul]\ar@{-}[dl]\\ & \bullet & }}\,,
\end{align}
where $\circ$ represents $\frac{d}{dx}$ and $\bullet$ represents $\frac{d}{dy}$.
Now, to compute $H^1_{P_0}(\Cxy)$, we have to consider elements
$f,g\in\Cxy$, satisfying the two independent equations: $\Phi_1(f)=0$ and $\Phi_2(g)=0$.
We first consider the second equation, $\Phi_2(g)=0$. We then have
\begin{prop}\label{prp:ker_phi_2}
The kernel of the linear map $\Phi_2$, from $\Cxy$ to $(T_{\Cxy}/[T_{\Cxy},T_{\Cxy}])_2$ is
$$
\ker(\Phi_2) = \CC[y].
$$
\end{prop}
\begin{proof}
Let $g\in\Cxy$ be polynomial in $x$ and $y$, such that $\Phi_2(g)=0$. 
Let us write $g = x g_0 + yg_1 +a$, where
$g_0, g_1 \in\Cxy$ and $a\in\CC$. Then, we have 
$\displaystyle\frac{dg}{dx} = 1\otimes g_0 + x\frac{dg_0}{dx} + y\frac{dg_1}{dx}$
and the equation $\Phi_2(g)=0$
becomes:
\begin{align*}
0 = \Phi_2(g) = & 
\,\raisebox{1.8pc}{\xymatrix@R=.75pc@C=.75pc{& \circ & \\ 
\lck{x}\ar@{-}[ur]\ar@{-}[dr] & & g_0
\ar@{-}[ul]\ar@{-}[dl]\\ & \bullet & }}
+\raisebox{1.8pc}{\xymatrix@R=.75pc@C=.75pc{& \circ & \\ 
\lck{x\frac{dg_0}{dx}'x}\ar@{-}[ur]\ar@{-}[dr] & & \frac{dg_0}{dx}''
\ar@{-}[ul]\ar@{-}[dl]\\ & \bullet & }}  
+\raisebox{1.8pc}{\xymatrix@R=.75pc@C=.75pc{& \circ & \\ 
\lck{y\frac{dg_1}{dx}'x}\ar@{-}[ur]\ar@{-}[dr] & & \frac{dg_1}{dx}''
\ar@{-}[ul]\ar@{-}[dl]\\ & \bullet & }}\\ &
-\raisebox{1.8pc}{\xymatrix@R=.75pc@C=.75pc{& \circ & \\ 
\lck{1}\ar@{-}[ur]\ar@{-}[dr] & & xg_0
\ar@{-}[ul]\ar@{-}[dl]\\ & \bullet & }}
- \raisebox{1.8pc}{\xymatrix@R=.75pc@C=.75pc{& \circ & \\ 
\lck{x\frac{dg_0}{dx}'}\ar@{-}[ur]\ar@{-}[dr] & & x\frac{dg_0}{dx}''
\ar@{-}[ul]\ar@{-}[dl]\\ & \bullet & }}  
- \raisebox{1.8pc}{\xymatrix@R=.75pc@C=.75pc{& \circ & \\ 
\lck{y\frac{dg_1}{dx}'}\ar@{-}[ur]\ar@{-}[dr] & & x\frac{dg_1}{dx}''
\ar@{-}[ul]\ar@{-}[dl]\\ & \bullet & }}\,.
\end{align*}
Then, we see that the term 
$\raisebox{1.8pc}{\xymatrix@R=.75pc@C=.75pc{& \circ & \\ 
\lck{1}\ar@{-}[ur]\ar@{-}[dr] & & xg_0
\ar@{-}[ul]\ar@{-}[dl]\\ & \bullet & }}$ has to cancel itself,
which means that $xg_0=0$ and $g = yg_1 + a\in y\,\Cxy +\CC$. 

Now, we know that we can write 
$g = a +\sum_{k\geq 1} y^k g_k$ with, for any $k\geq 1$, $g_k\in x\,\Cxy+ \CC$. We then have 
$\displaystyle\frac{dg}{dx} = \sum_{k\geq 1} y^k\, \frac{dg_k}{dx}
= \sum_{k\geq 1} y^k \,\frac{dg_k}{dx}'\otimes\frac{dg_k}{dx}''$ and
the equation $\Phi_2(g)=0$ becomes
$$
\Phi_2(g) = 0 = \sum_{k\geq 1} y^k\,\Phi_2(g_k).
$$
So that, for each $k\geq 1$, we must have $\Phi_2(g_k)=0$. But we have
seen above that this implies $g_k\in y\,\Cxy + \CC$, while we have
assumed that $g_k \in x\,\Cxy +\CC$, thus $g_k\in\CC$, for all
$k\in\N^*$. We then conclude that $g\in\CC[y]$.
\end{proof}

Now, let us study the first equation $\Phi_1(f)=0$. 
To do this, we will need the following
\begin{lemma}\label{lma:phi_1(m,n)}
Let $s\in\N^*$ and $(k_1,k_2,k_3,\dots,k_s)\in(\N^*)^{2s}$. Let 
$m = x^{k_1}y^{k_2}x^{k_3}\cdots x^{k_{2s-1}}y^{k_{2s}} \in \Cxy$. 
We have
\begin{align}\label{eqn:phi_1(m)}
\Phi_1(m) = & 
 \sum_{i=1}^s\; 
\,\raisebox{1.8pc}{\xymatrix@R=.75pc@C=.75pc{& \circ & \\ 
\lck{x^{k_1}\cdots y^{k_{2(i-1)}}x^{k_{(2i-1)}}}\ar@{-}[ur]\ar@{-}[dr] & & 
   y^{k_{2i}}\cdots y^{k_{2s}}
\ar@{-}[ul]\ar@{-}[dl]\\ & \circ & }}\nonumber\\ &
-\sum_{i=2}^s\; 
\,\raisebox{1.8pc}{\xymatrix@R=.75pc@C=.75pc{& \circ & \\ 
\lck{x^{k_1}\cdots y^{k_{2(i-1)}}}\ar@{-}[ur]\ar@{-}[dr] & & 
   x^{k_{(2i-1)}}y^{k_{2i}}\cdots y^{k_{2s}}
\ar@{-}[ul]\ar@{-}[dl]\\ & \circ & }}\,.
\end{align}
Next, let $n = y^{k_1}x^{k_2}y^{k_3}\cdots
y^{k_{2s-1}}x^{k_{2s}} \in \Cxy$, then
\begin{align}\label{eqn:phi_1(n)}
\Phi_1(n) = & 
- \sum_{i=1}^s\; 
\,\raisebox{1.8pc}{\xymatrix@R=.75pc@C=.75pc{& \circ & \\ 
\lck{y^{k_1}\cdots x^{k_{2(i-1)}}y^{k_{(2i-1)}}}\ar@{-}[ur]\ar@{-}[dr] & & 
   x^{k_{2i}}\cdots x^{k_{2s}}
\ar@{-}[ul]\ar@{-}[dl]\\ & \circ & }}\nonumber\\ &
+\sum_{i=2}^s\; 
\,\raisebox{1.8pc}{\xymatrix@R=.75pc@C=.75pc{& \circ & \\ 
\lck{y^{k_1}\cdots x^{k_{2(i-1)}}}\ar@{-}[ur]\ar@{-}[dr] & & 
   y^{k_{(2i-1)}}x^{k_{2i}}\cdots x^{k_{2s}}
\ar@{-}[ul]\ar@{-}[dl]\\ & \circ & }}\,.
\end{align}
\end{lemma}
\begin{proof}
We will prove the first statement of the lemma. The proof of the second statement is completely analogous.
Let $m=x^{k_1}y^{k_2}x^{k_3}\cdots x^{k_{2s-1}}y^{k_{2s}}\in\Cxy$, 
where $s\in\N^*$ and where $k_1,k_2,k_3,\dots,k_s\in\N^*$ are (non zero) integers.
First of all, we compute:
\begin{equation*}
\frac{dm}{dx} = 
\sum_{i=1}^s\; \sum_{j=0}^{k_{(2i-1)}-1} 
  x^{k_1}\cdots y^{k_{2(i-1)}}x^j\otimes
  x^{(k_{(2i-1)}-1-j)}y^{k_{2i}}\cdots y^{k_{2s}}.
\end{equation*}
Then, by definition of $\Phi_1$, we have
\begin{align*}
\Phi_1(m) = & 
-\,\raisebox{1.8pc}{\xymatrix@R=.75pc@C=.75pc{& \circ & \\ 
\lck{x^{k_1}\dots y^{k_{2s}}}\ar@{-}[ur]\ar@{-}[dr] & & 1
    \ar@{-}[ul]\ar@{-}[dl]\\ & \circ & }} \\ &
+\sum_{i=1}^s\; \sum_{j=0}^{k_{(2i-1)}-1}
\,\raisebox{1.8pc}{\xymatrix@R=.75pc@C=.75pc{& \circ & \\ 
\lck{x^{k_1}\cdots y^{k_{2(i-1)}}x^{j+1}}\ar@{-}[ur]\ar@{-}[dr] & & 
   x^{(k_{(2i-1)}-1-j)}y^{k_{2i}}\cdots y^{k_{2s}}
\ar@{-}[ul]\ar@{-}[dl]\\ & \circ & }}\\ &
-\sum_{i=1}^s\; \sum_{j=0}^{k_{(2i-1)}-1}
\,\raisebox{1.8pc}{\xymatrix@R=.75pc@C=.75pc{& \circ & \\ 
\lck{x^{k_1}\cdots y^{k_{2(i-1)}}x^j}\ar@{-}[ur]\ar@{-}[dr] & & 
   x^{(k_{(2i-1)}-j)}y^{k_{2i}}\cdots y^{k_{2s}}
\ar@{-}[ul]\ar@{-}[dl]\\ & \circ & }}\\ 
=& -\,\raisebox{1.8pc}{\xymatrix@R=.75pc@C=.75pc{& \circ & \\ 
\lck{x^{k_1}\dots y^{k_{2s}}}\ar@{-}[ur]\ar@{-}[dr] & & 1
    \ar@{-}[ul]\ar@{-}[dl]\\ & \circ & }} \\ &
+ \sum_{i=1}^s\; 
\,\raisebox{1.8pc}{\xymatrix@R=.75pc@C=.75pc{& \circ & \\ 
\lck{x^{k_1}\cdots y^{k_{2(i-1)}}x^{k_{(2i-1)}}}\ar@{-}[ur]\ar@{-}[dr] & & 
   y^{k_{2i}}\cdots y^{k_{2s}}
\ar@{-}[ul]\ar@{-}[dl]\\ & \circ & }}\\ &
-\sum_{i=1}^s\; 
\,\raisebox{1.8pc}{\xymatrix@R=.75pc@C=.75pc{& \circ & \\ 
\lck{x^{k_1}\cdots y^{k_{2(i-1)}}}\ar@{-}[ur]\ar@{-}[dr] & & 
   x^{k_{(2i-1)}}y^{k_{2i}}\cdots y^{k_{2s}}
\ar@{-}[ul]\ar@{-}[dl]\\ & \circ & }}\,,
\end{align*}
which yields the formula (\ref{eqn:phi_1(m)}). 
\end{proof}
We will also need the following formula that gives a nice
interpretation of $d_{P_0}(m)$, where $m$ is a monomial like in the
previous lemma.
\begin{lemma}\label{lma:d^0(m)}
Let $s\in\N^*$ and $(k_1,k_2,k_3,\dots,k_s)\in(\N^*)^{2s}$. We 
consider the monomial in $\Cxy$, written as
$m = x^{k_1}y^{k_2}x^{k_3}\cdots x^{k_{2s-1}}y^{k_{2s}}$. 
Then, we have 
\begin{eqnarray*}\label{eqn:d^0_P0(m)}
\frac{dm}{dx}''\frac{dm}{dx}'x - x \frac{dm}{dx}''\frac{dm}{dx}' &=&
- \sum_{i=1}^{s} x^{k_{2i-1}}\cdots
  x^{k_{2s-1}}y^{k_{2s}}x^{k_1}y^{k_2}\cdots y^{k_{2(i-1)}}\nonumber\\
&&+\sum_{i=1}^s\; 
y^{k_{2i}}\cdots y^{k_{2s}}x^{k_1}\cdots y^{k_{2(i-1)}}x^{k_{(2i-1)}}.
\end{eqnarray*}
That is, $d_{P_0}^0(m)$ is obtained by considering all the
cyclic permutations of the blocks $x^j$ and $y^j$ in $m$ (together with the sign of
the permutation) and multiplying the result by $\frac{d}{dx}$.
\end{lemma}
\begin{proof}
Similar to the proof of the previous lemma, we have
\begin{equation*}
\frac{dm}{dx} = 
\sum_{i=1}^s\; \sum_{j=0}^{k_{(2i-1)}-1} 
x^{k_1}\cdots y^{k_{2(i-1)}}x^j\otimes x^{(k_{(2i-1)}-1-j)}y^{k_{2i}}\cdots y^{k_{2s}}.
\end{equation*}
From this, it is straightforward to obtain equation (\ref{eqn:d^0_P0(m)}).
\end{proof}

We are now able to determine the first double Poisson cohomology group of
the double Poisson algebra $(\Cxy, P_0)$.
\begin{prop}
Let us consider the linear double Poisson tensor $P_0 = x
\frac{d}{dx}\frac{d}{dx}$ on $\Cxy$. The first double Poisson cohomology space,
associated to $P_0$ is given by:
\begin{equation*}
H^1_{P_0} \simeq \CC\,\frac{d}{dx} \,\oplus\, \CC[y]\frac{d}{dy}.
\end{equation*}
\end{prop}
\begin{proof}
Let $f\frac{d}{dx}+g\frac{d}{dy}$ be a $1$-cocycle of the double
Poisson cohomology associated to the double Poisson algebra
$(\Cxy,P_0)$. We have seen that the cocycle condition can be written
as:
\begin{eqnarray*}
\left\lbrace
 \begin{array}{ccl}
  \Phi_1(f) = 0,\\
  \Phi_2(g) = 0,
 \end{array}
\right.
\end{eqnarray*}
where the operators $\Phi_1$ and $\Phi_2$ are defined in
(\ref{def:Phi_1}) and (\ref{def:Phi_2}). According to Proposition
\ref{prp:ker_phi_2}, we know that $g\in\CC[y]$. As for any $h\in\Cxy$,
$d_{P_0}(h)\in\Cxy\frac{d}{dx}$, it is clear that
the elements of $\CC[y]\frac{d}{dy}$ give non-trivial double Poisson
cohomology classes in $H^1_{P_0}(\Cxy)$. 

Remains to consider the equation $\Phi_1(f) = 0$. First of all observe
this equation implies $f\in x\Cxy y + y\Cxy x + \CC$. In fact, 
suppose that there is a monomial in $f$ that can be written as
$xf_0\,x$, where $f_0\in\Cxy$. Then, we have
\begin{align*}
\Phi_1(xf_0\,x) = & 
 \,\raisebox{1.8pc}{\xymatrix@R=.75pc@C=.75pc{& \circ & \\ 
\lck{x}\ar@{-}[ur]\ar@{-}[dr] & & f_0 x
\ar@{-}[ul]\ar@{-}[dl]\\ & \circ & }}
+\,\raisebox{1.8pc}{\xymatrix@R=.75pc@C=.75pc{& \circ & \\ 
\lck{x\frac{df_0}{dx}'x}\ar@{-}[ur]\ar@{-}[dr] & & \frac{df_0}{dx}''x
\ar@{-}[ul]\ar@{-}[dl]\\ & \circ & }} 
+ \,\raisebox{1.8pc}{\xymatrix@R=.75pc@C=.75pc{& \circ & \\ 
\lck{xf_0\, x}\ar@{-}[ur]\ar@{-}[dr] & & 1
\ar@{-}[ul]\ar@{-}[dl]\\ & \circ & }} \\ &
- \,\raisebox{1.8pc}{\xymatrix@R=.75pc@C=.75pc{& \circ & \\ 
\lck{x\frac{df_0}{dx}'}\ar@{-}[ur]\ar@{-}[dr] & & x\frac{df_0}{dx}''x
\ar@{-}[ul]\ar@{-}[dl]\\ & \circ & }}
-\,\raisebox{1.8pc}{\xymatrix@R=.75pc@C=.75pc{& \circ & \\ 
\lck{xf_0}\ar@{-}[ur]\ar@{-}[dr] & & x
\ar@{-}[ul]\ar@{-}[dl]\\ & \circ & }}.
\end{align*}
But then $\Phi_1(f) = 0$, implies that the term 
$\raisebox{1.8pc}{\xymatrix@R=.75pc@C=.75pc{& \circ & \\ 
\lck{xf_0\, x}\ar@{-}[ur]\ar@{-}[dr] & & 1
\ar@{-}[ul]\ar@{-}[dl]\\ & \circ & }} $ has to cancel itself, so that
$xf_0\,x$ has to be zero. Now suppose a monomial of the form $yf_0\,y$
appears in $f$. We have
\begin{align*}
\Phi_1(y f_0\,y) = & 
-\,\raisebox{1.8pc}{\xymatrix@R=.75pc@C=.75pc{& \circ & \\ 
\lck{yf_0\,y}\ar@{-}[ur]\ar@{-}[dr] & & 1\ar@{-}[ul]\ar@{-}[dl]\\ & \circ & }}
+\,\raisebox{1.8pc}{\xymatrix@R=.75pc@C=.75pc{& \circ & \\ 
\lck{y\frac{df_0}{dx}'x}\ar@{-}[ur]\ar@{-}[dr] & & \frac{df_0}{dx}''y
\ar@{-}[ul]\ar@{-}[dl]\\ & \circ & }}
-\,\raisebox{1.8pc}{\xymatrix@R=.75pc@C=.75pc{& \circ & \\ 
\lck{y\frac{df_0}{dx}'}\ar@{-}[ur]\ar@{-}[dr] & & x\frac{df_0}{dx}''y
\ar@{-}[ul]\ar@{-}[dl]\\ & \circ & }}\,.
\end{align*}
The term $\raisebox{1.8pc}{\xymatrix@R=.75pc@C=.75pc{& \circ & \\ 
\lck{yf_0\,y}\ar@{-}[ur]\ar@{-}[dr] & & 1\ar@{-}[ul]\ar@{-}[dl]\\ &
\circ & }}$ cannot appear in $\Phi_1(f)$ in any other way and hence
has to vanish, i.e. $yf_0\,y$ has to be zero. 

But then we know that $f$ can be written as $f = \sum\limits_{s\in\N^*}
f_{2s} + a$,
where $a\in\CC$ and $f_{2s} \in\Cxy$ is of the form
$$
f_{2s} := \sum_{\stackrel{K=(k_1,\dots,k_{2s})}{\in(\N^*)^{2s}}} 
c_K\; x^{k_1}y^{k_2}x^{k_3}\cdots x^{k_{2s-1}}y^{k_{2s}}
-
\sum_{\stackrel{L=(\ell_1,\dots,\ell_{2s})}{\in(\N^*)^{2s}}} 
\tilde c_L\; y^{\ell_1}x^{\ell_2}y^{\ell_3}\cdots y^{\ell_{2s-1}}x^{\ell_{2s}},
$$
where $c_K$ and $\tilde c_L$ are constants. According to Lemma
\ref{lma:phi_1(m,n)}, the equation $\Phi_1(f)=0$ implies that, for
each $s\in\N^*$, $\Phi_1(f_{2s})=0$ (i.e., $\Phi_1(f_{2s})$ can not be
cancelled by another $\Phi_1(f_{2s'})$).  

Let us then consider the equation $\Phi_1(f_{2s})=0$.
According to Lemma \ref{lma:phi_1(m,n)}, by collecting the terms of
the form $\raisebox{1.8pc}{\xymatrix@R=.75pc@C=.75pc{& \circ & \\ 
\lck{x\cdots x}\ar@{-}[ur]\ar@{-}[dr] & & y\cdots y\ar@{-}[ul]\ar@{-}[dl]\\ &
\circ & }}$ and of the form $\raisebox{1.8pc}{\xymatrix@R=.75pc@C=.75pc{& \circ & \\ 
\lck{x\cdots y}\ar@{-}[ur]\ar@{-}[dr] & & x\cdots y\ar@{-}[ul]\ar@{-}[dl]\\ &
\circ & }}$ (which have to be cancelled by terms of the same form), 
we get the three following equations:
\begin{align}\label{eqn:phi_1_1}
0 = &\sum_{\stackrel{K=(k_1,\dots,k_{2s})}{\in(\N^*)^{2s}}} 
c_K\;\sum_{i=1}^s\; 
\,\raisebox{1.8pc}{\xymatrix@R=.75pc@C=.75pc{& \circ & \\ 
\lck{x^{k_1}\cdots y^{k_{2(i-1)}}x^{k_{(2i-1)}}}\ar@{-}[ur]\ar@{-}[dr] & & 
   y^{k_{2i}}\cdots y^{k_{2s}}
\ar@{-}[ul]\ar@{-}[dl]\\ & \circ & }}\nonumber\\ 
+ &\sum_{\stackrel{L=(\ell_1,\dots,\ell_{2s})}{\in(\N^*)^{2s}}}\tilde c_L\;
\sum_{i=1}^s\; 
\,\raisebox{1.8pc}{\xymatrix@R=.75pc@C=.75pc{& \circ & \\ 
\lck{y^{k_1}\cdots x^{k_{2(i-1)}}y^{k_{(2i-1)}}}\ar@{-}[ur]\ar@{-}[dr] & & 
   x^{k_{2i}}\cdots x^{k_{2s}}
\ar@{-}[ul]\ar@{-}[dl]\\ & \circ & }}
\end{align}
and
\begin{align}\label{eqn:phi_1_2}
0 = &\sum_{\stackrel{K=(k_1,\dots,k_{2s})}{\in(\N^*)^{2s}}} 
c_K\; \sum_{i=2}^s\; 
\,\raisebox{1.8pc}{\xymatrix@R=.75pc@C=.75pc{& \circ & \\ 
\lck{x^{k_1}\cdots y^{k_{2(i-1)}}}\ar@{-}[ur]\ar@{-}[dr] & & 
   x^{k_{(2i-1)}}y^{k_{2i}}\cdots y^{k_{2s}}
\ar@{-}[ul]\ar@{-}[dl]\\ & \circ & }}\,.
\end{align}

From Equation (\ref{eqn:phi_1_1}), we conclude, as the first sum cannot cancel itself, that, for each $1\leq i\leq s$
\begin{align*}
&\sum_{\stackrel{K=(k_1,\dots,k_{2s})}{\in(\N^*)^{2s}}} 
c_K\; 
\,\raisebox{1.8pc}{\xymatrix@R=.75pc@C=.75pc{& \circ & \\ 
\lck{x^{k_1}\cdots y^{k_{2(i-1)}}x^{k_{(2i-1)}}}\ar@{-}[ur]\ar@{-}[dr] & & 
   y^{k_{2i}}\cdots y^{k_{2s}}
\ar@{-}[ul]\ar@{-}[dl]\\ & \circ & }}\\
= & -\sum_{\stackrel{K=(k_1,\dots,k_{2s})}{\in(\N^*)^{2s}}} 
c_K\; 
\,\raisebox{1.8pc}{\xymatrix@R=.75pc@C=.75pc{& \circ & \\ 
\lck{y^{k_{2i}}\cdots y^{k_{2s}}}\ar@{-}[ur]\ar@{-}[dr] & & 
   x^{k_1}\cdots x^{k_{(2i-1)}}
\ar@{-}[ul]\ar@{-}[dl]\\ & \circ & }}\\
= & -\sum_{\stackrel{L=(\ell_1,\dots,\ell_{2s})}{\in(\N^*)^{2s}}}\tilde c_L\;
\,\raisebox{1.8pc}{\xymatrix@R=.75pc@C=.75pc{& \circ & \\ 
\lck{y^{k_1}\cdots x^{k_{2(i-1)}}y^{k_{(2i-1)}}}\ar@{-}[ur]\ar@{-}[dr] & & 
   x^{k_{2i}}\cdots x^{k_{2s}}
\ar@{-}[ul]\ar@{-}[dl]\\ & \circ & }}
\end{align*}
and this can only happen if, for each $1\leq i\leq s$: 
\begin{eqnarray}\label{eqn:phi_11}
\sum_{\stackrel{K=(k_1,\dots,k_{2s})}{\in(\N^*)^{2s}}} 
c_K\; y^{k_{2i}}\cdots y^{k_{2s}}  x^{k_1}\cdots x^{k_{(2i-1)}}
= \sum_{\stackrel{L=(\ell_1,\dots,\ell_{2s})}{\in(\N^*)^{2s}}}\tilde c_L\;
 y^{k_1}\cdots x^{k_{2s}}.
\end{eqnarray}
Then, in the equation (\ref{eqn:phi_1_2}), the sum obtained for 
$2\leq i\leq s$ has to be cancelled by the one obtained for the $s-i+2$,
i.e.,
\begin{align*}
&\sum_{\stackrel{K=(k_1,\dots,k_{2s})}{\in(\N^*)^{2s}}} 
c_K\; 
\,\raisebox{1.8pc}{\xymatrix@R=.75pc@C=.75pc{& \circ & \\ 
\lck{x^{k_1}\cdots y^{k_{2(i-1)}}}\ar@{-}[ur]\ar@{-}[dr] & & 
   x^{k_{(2i-1)}}y^{k_{2i}}\cdots y^{k_{2s}}
\ar@{-}[ul]\ar@{-}[dl]\\ & \circ & }}\\
= &- \sum_{\stackrel{K=(k_1,\dots,k_{2s})}{\in(\N^*)^{2s}}} 
c_K\; 
\,\raisebox{1.8pc}{\xymatrix@R=.75pc@C=.75pc{& \circ & \\ 
\lck{x^{k_{(2i-1)}}y^{k_{2i}}\cdots y^{k_{2s}}}\ar@{-}[ur]\ar@{-}[dr] & & 
   x^{k_1}\cdots y^{k_{2(i-1)}}
\ar@{-}[ul]\ar@{-}[dl]\\ & \circ & }}\\
= & \sum_{\stackrel{K=(k_1,\dots,k_{2s})}{\in(\N^*)^{2s}}} 
c_K\; 
\,\raisebox{1.8pc}{\xymatrix@R=.75pc@C=.75pc{& \circ & \\ 
\lck{x^{k_1}\cdots y^{k_{2(s-i+1)}}}\ar@{-}[ur]\ar@{-}[dr] & & 
   x^{k_{(2(s-i)+3)}}y^{k_{2(s-i+2)}}\cdots y^{k_{2s}}
\ar@{-}[ul]\ar@{-}[dl]\\ & \circ & }}
\end{align*}
when written with exactly $2(s-i+1)$ blocks of the form $x^j$ or $y^j$ in
the box. This implies, for each $2\leq i \leq s$, that
\begin{eqnarray}\label{eqn:phi_12}
-\sum_{\stackrel{K=(k_1,\dots,k_{2s})}{\in(\N^*)^{2s}}} 
c_K\; x^{k_{(2i-1)}}\cdots y^{k_{2s}}x^{k_1}\cdots y^{k_{2(i-1)}}
 = \sum_{\stackrel{K=(k_1,\dots,k_{2s})}{\in(\N^*)^{2s}}} 
c_K\; x^{k_1}\cdots y^{k_{2s}}.
\end{eqnarray}
Now let
$$
h_{2s} := \sum_{\stackrel{K=(k_1,\dots,k_{2s})}{\in(\N^*)^{2s}}}
c_K\; x^{k_1}\cdots y^{k_{2s}} \in\Cxy.
$$
According to Lemma \ref{lma:d^0(m)}, we have:
\begin{eqnarray*}
\lefteqn{\frac{dh_{2s}}{dx}''\frac{dh_{2s}}{dx}'x 
 - x \frac{dh_{2s}}{dx}''\frac{dh_{2s}}{dx}' =}\\
&&- \sum_{i=1}^{s}\,\sum_{\stackrel{K=(k_1,\dots,k_{2s})}{\in(\N^*)^{2s}}} 
c_K\; x^{k_{2i-1}}\cdots
  x^{k_{2s-1}}y^{k_{2s}}x^{k_1}y^{k_2}\cdots y^{k_{2(i-1)}}\\
&&+\sum_{i=1}^s\,\sum_{\stackrel{K=(k_1,\dots,k_{2s})}{\in(\N^*)^{2s}}} 
c_K\; 
y^{k_{2i}}\cdots y^{k_{2s}}x^{k_1}\cdots y^{k_{2(i-1)}}x^{k_{(2i-1)}}.
\end{eqnarray*}
From equations (\ref{eqn:phi_11}) and (\ref{eqn:phi_12}), we obtain
$$
\frac{dh_{2s}}{dx}''\frac{dh_{2s}}{dx}'x - x \frac{dh_{2s}}{dx}''\frac{dh_{2s}}{dx}' =
-s\, f_{2s}.
$$
According to Proposition \ref{imageofdp00}, this yields
$$
f\frac{d}{dx}=\sum_{s\in\N^*} f_{2s}\frac{d}{dx} + a\frac{d}{dx} 
 = d_{P_0}\left(-\frac{1}{s}\,h_{2s}\right) +a\frac{d}{dx},
$$
and we conclude that
$
H^1_{P_0}(\Cxy) \simeq \CC\,\frac{d}{dx}\oplus\CC[y]\, \frac{d}{dy}.$
\end{proof}

An analogous proof shows that
\begin{prop}
Let us consider the linear double Poisson tensor 
$$
\tilde P_0 :=P_{lin}^{\CC\times\CC}= x
\frac{d}{dx}\frac{d}{dx} + y\frac{d}{dy}\frac{d}{dy}
$$
on $\Cxy$. 
Then, we have:
\begin{equation*}
H_{\tilde P_0}^0(\Cxy) \simeq \CC[x]\oplus\CC[y]
\end{equation*}
and the first double Poisson cohomology space,
associated to $\tilde P_0$ is given by:
\begin{equation*}
H^1_{\tilde P_0}(\Cxy) \simeq \CC\,\frac{d}{dx} \,\oplus\, \CC\,\frac{d}{dy}.
\end{equation*}
\end{prop}
\begin{remark}
On $\rep_1(\Cxy)$, the double Poisson tensors $P_0=P_{lin}^{\CC\times\CC\epsilon^2}$ and $\tilde P_0=P_{lin}^{\CC\times\CC}$ induce the trivial Poisson bracket.
However, on $\rep_n(\Cxy)$ ($n\geq 2$) the double Poisson tensor $P_{lin}^{\CC\times\CC}$
is mapped (by the trace map) to the canonical linear Poisson structure on the product $\gl_n^*\times \gl_n^*$. For the Lie algebra $\gl_n$, we know
\cite{Fuks} that the Lie algebra cohomology space $H^1_L(\gl_n;\CC)$ is of dimension $1$. To obtain the first Poisson cohomology group of $\gl_n^*\times
\gl_n^*$, we have to consider the tensor product of $H^1_L(\gl_n\times \gl_n;\CC)$ which is of dimension two and the algebra of the Casimirs
of $\gl_n^*\times \gl_n^*$, which is an infinite dimensional vector space. That is, the trace map from $H^1_{\tilde P_0}(\Cxy)$ to $H^1_{tr(\tilde P_0)}(\rep_n(\Cxy))$ is not onto.
\end{remark}
\subsection{The linear double Poisson tensors 
$x \frac{d}{dx}\frac{d}{dx} + y\frac{d}{dx}\frac{d}{dy}$, 
$x \frac{d}{dx}\frac{d}{dy} + y\frac{d}{dy}\frac{d}{dy}$}

Now we consider the linear double Poisson tensor: 
$$
P_1:=P_{lin}^{B_2^1}= x \frac{d}{dx}\frac{d}{dx} +
y\frac{d}{dx}\frac{d}{dy}.
$$
 We will determine the double Poisson
cohomology groups $H^0_{P_1}(\Cxy)$ and
$H^1_{P_1}(\Cxy)$.
We begin by observing
\begin{lemma}\label{lma:gauge}
Let us consider the free algebra $\CC\langle x_1,\dots,x_n\rangle$,
associated to the quiver $Q$, with one vertex and $n$ loops $x_1,\dots,x_n$. For each
$h\in \CC\langle x_1,\dots,x_n\rangle$, we have  
\begin{eqnarray*}
\sum_{i=1}^n \left(\left(\frac{d}{dx_i}\circ x_i\right)(h)
-\left(x_i\circ\frac{d}{dx_i}\right)(h)\right)= h\otimes 1 -1\otimes h,
\end{eqnarray*}
(where $\circ$ denotes the inner multiplication). This can also be written as:
\begin{eqnarray*}
\sum_{i=1}^n \left(\left(\frac{dh}{dx_i}\right)'x_i\otimes \left(\frac{dh}{dx_i}\right)''
-\left(\frac{dh}{dx_i}\right)'\otimes x_i\left(\frac{dh}{dx_i}\right)''\right)
= h\otimes 1 -1\otimes h.
\end{eqnarray*}
\end{lemma}
\begin{proof}
This can easily be seen from the definition of the $\frac{d}{dx_i}$,
but it is also a particular case of Proposition 6.2.2 of
\cite{MichelDPA}, which states that the gauge element $E$ of $Q$ is
given by $\displaystyle E=\sum_{a\in Q_1} \left\lbrack \frac{d}{da},a\right\rbrack$.
\end{proof}

Now, let us first consider the double Poisson cohomology space
$H^0_{P_1}(\Cxy)$.
\begin{prop}\label{imageofdp10}
For $f\in\Cxy$, we have
$$d_{P_1}^0(f) = 
\xymatrix@R=1.75pc{\lck{\circ}\ar@{-}@/^.25pc/[r] & y\frac{df}{dy}''\frac{df}{dy}'\ar@{-}@/^.5pc/[l]}
-
\xymatrix@R=.75pc{\lck{\bullet}\ar@{-}@/^.25pc/[r] & y\frac{df}{dx}''\frac{df}{dx}'\ar@{-}@/^.5pc/[l]}.
$$
Which means that
$$H_{P_1}^0(\Cxy) = \CC.$$
\end{prop}
\begin{proof}
In fact, by computing $d_{P_1}(f) = \{P_1,f\}$, one obtains exactly:
$$
d_{P_1}^0(f) = 
\xymatrix@R=1.75pc{\lck{\circ}\ar@{-}@/^.25pc/[r] & 
\left(\frac{df}{dx}''\frac{df}{dx}'x - x\frac{df}{dx}''\frac{df}{dx}'+\frac{df}{dy}''\frac{df}{dy}'y\right)\ar@{-}@/^.5pc/[l]}
-
\xymatrix@R=.75pc{\lck{\bullet}\ar@{-}@/^.25pc/[r] & y\frac{df}{dx}''\frac{df}{dx}'\ar@{-}@/^.5pc/[l]}.
$$
According to Lemma \ref{lma:gauge}, we have 
$$
\frac{df}{dx}' x\otimes\frac{df}{dx}'' +\frac{df}{dy}'y\otimes
\frac{df}{dy}''+ 1\otimes f
= \frac{df}{dx}'\otimes x \frac{df}{dx}'' + \frac{df}{dy}'\otimes
y\frac{df}{dy}'' + f\otimes 1.
$$
Applying $-^{op}$ and $\mu$, this last formula gives:
$$
\frac{df}{dx}''\frac{df}{dx}'x +\frac{df}{dy}''\frac{df}{dy}'y 
= x\frac{df}{dx}''\frac{df}{dx}' + y\frac{df}{dy}''\frac{df}{dy}',
$$
which leads to the expression for $d_{P_1}(f)$ given above.

Suppose now that $f$ is a $0$-cocycle, that is to say
$d_{P_1}(f)=0$. This is equivalent to say that 
$$
y\frac{df}{dy}''\frac{df}{dy}' = y\frac{df}{dx}''\frac{df}{dx}' = 0,
$$
which yields
$$
\frac{df}{dy}''\frac{df}{dy}' = \frac{df}{dx}''\frac{df}{dx}' = 0
$$
and 
$$
y\frac{df}{dy}''\frac{df}{dy}' + x \frac{df}{dx}''\frac{df}{dx}' = 0.
$$
This implies
$$
\frac{df}{dy}'y\frac{df}{dy}'' + \frac{df}{dx}'x \frac{df}{dx}'' \in \lbrack\Cxy,\Cxy\rbrack,
$$
Using the NC-Euler Formula (Proposition \ref{nceuler}),
we can then write 
$$f\in \CC \oplus [\Cxy,\Cxy].$$
Finally, $H^0_{P_1}(\Cxy)= \ker(d_{P_1}^0)/[\Cxy,\Cxy]= \CC$, which concludes the proof. 
\end{proof}
Let us now determine $H^1_{P_1}(\Cxy)$. We will first use Lemma \ref{lma:gauge} to obtain a useful expression for the coboundary
operator $d_{P_1}^1$.
\begin{lemma}\label{imageofdp11}
Let $\displaystyle f\frac{d}{dx}+g\frac{d}{dy}\in (T_{\Cxy}/[T_{\Cxy},T_{\Cxy}])_1$, then
\begin{align*}
d_{P_1}^1\left(f\frac{d}{dx}+g\frac{d}{dy}\right) = & 
-\,\raisebox{1.8pc}{\xymatrix@R=.75pc@C=.75pc{& \circ & \\ 
\lck{1}\ar@{-}[ur]\ar@{-}[dr] & & f\ar@{-}[ul]\ar@{-}[dl]\\ & \circ & }}
+ \raisebox{1.8pc}{\xymatrix@R=.75pc@C=.75pc{& \circ & \\ \lck{\frac{df}{dy}'}\ar@{-}[ur]\ar@{-}[dr] & & y\frac{df}{dy}''\ar@{-}[ul]\ar@{-}[dl]\\ & \circ & }}
- \raisebox{1.8pc}{\xymatrix@R=.75pc@C=.75pc{& \circ & \\ \lck{1}\ar@{-}[ur]\ar@{-}[dr] & & g\ar@{-}[ul]\ar@{-}[dl]\\ & \bullet & }}  \\
&
+ \raisebox{1.8pc}{\xymatrix@R=.75pc@C=.75pc{& \circ & \\ \lck{\frac{dg}{dy}'}\ar@{-}[ur]\ar@{-}[dr] & & y\frac{dg}{dy}''\ar@{-}[ul]\ar@{-}[dl]\\ & \bullet & }}
+ \raisebox{1.8pc}{\xymatrix@R=.75pc@C=.75pc{& \circ & \\ \lck{y\frac{df}{dx}''}\ar@{-}[ur]\ar@{-}[dr] & & \frac{df}{dx}'\ar@{-}[ul]\ar@{-}[dl]\\ & \bullet & }}
- \raisebox{1.8pc}{\xymatrix@R=.75pc@C=.75pc{& \bullet& \\ \lck{\frac{dg}{dx}'}\ar@{-}[ur]\ar@{-}[dr] & & y\frac{dg}{dx}''\,.\ar@{-}[ul]\ar@{-}[dl]\\ & \bullet & }}
\end{align*}
\end{lemma}
\begin{proof}
First, by computing 
$d_{P_1}\left(f\frac{d}{dx}+g\frac{d}{dy}\right)=\{P_1,f\frac{d}{dx}+g\frac{d}{dy}\}$,
one can write
\begin{equation*}
d_{P_1}\left(f\frac{d}{dx}+g\frac{d}{dy}\right) = (A) + (B) + (C),
\end{equation*}
where 
\begin{align*}
(A) = & 
\,-\raisebox{1.8pc}{\xymatrix@R=.75pc@C=.75pc{& \circ & \\ 
\lck{f}\ar@{-}[ur]\ar@{-}[dr] & & 1\ar@{-}[ul]\ar@{-}[dl]\\ & \circ & }}
+ \raisebox{1.8pc}{\xymatrix@R=.75pc@C=.75pc{& \circ & \\ \lck{\frac{df}{dx}'x}\ar@{-}[ur]\ar@{-}[dr] & & \frac{df}{dx}''\ar@{-}[ul]\ar@{-}[dl]\\ & \circ & }}
- \raisebox{1.8pc}{\xymatrix@R=.75pc@C=.75pc{& \circ & \\ \lck{\frac{df}{dx}'}\ar@{-}[ur]\ar@{-}[dr] & & x\frac{df}{dx}''\ar@{-}[ul]\ar@{-}[dl]\\ & \circ & }}
+ \raisebox{1.8pc}{\xymatrix@R=.75pc@C=.75pc{& \circ & \\
    \lck{\frac{df}{dy}'y}\ar@{-}[ur]\ar@{-}[dr] & & \frac{df}{dy}''\, ,\ar@{-}[ul]\ar@{-}[dl]\\ & \circ & }}
\end{align*}
\begin{align*}
(B)= &- \raisebox{1.8pc}{\xymatrix@R=.75pc@C=.75pc{& \circ & \\ \lck{g}\ar@{-}[ur]\ar@{-}[dr] & & 1\ar@{-}[ul]\ar@{-}[dl]\\ & \bullet & }}
+ \raisebox{1.8pc}{\xymatrix@R=.75pc@C=.75pc{& \circ & \\ \lck{\frac{dg}{dx}'x}\ar@{-}[ur]\ar@{-}[dr] & & \frac{dg}{dx}''\ar@{-}[ul]\ar@{-}[dl]\\ & \bullet & }}
- \raisebox{1.8pc}{\xymatrix@R=.75pc@C=.75pc{& \circ & \\
    \lck{\frac{dg}{dx}'}\ar@{-}[ur]\ar@{-}[dr] & &
    x\frac{dg}{dx}''\ar@{-}[ul]\ar@{-}[dl]\\ & \bullet & }}\\
&
+ \raisebox{1.8pc}{\xymatrix@R=.75pc@C=.75pc{& \circ & \\ \lck{\frac{dg}{dy}'y}\ar@{-}[ur]\ar@{-}[dr] & & \frac{dg}{dy}''\ar@{-}[ul]\ar@{-}[dl]\\ & \bullet & }}
+ \raisebox{1.8pc}{\xymatrix@R=.75pc@C=.75pc{& \circ & \\ \lck{y\frac{df}{dx}''}\ar@{-}[ur]\ar@{-}[dr] & & \frac{df}{dx}'\ar@{-}[ul]\ar@{-}[dl]\\ & \bullet & }}
\end{align*}
and
\begin{align*}
(C)=& - \raisebox{1.8pc}{\xymatrix@R=.75pc@C=.75pc{& \bullet& \\ \lck{\frac{dg}{dx}'}\ar@{-}[ur]\ar@{-}[dr] & & y\frac{dg}{dx}''\,.\ar@{-}[ul]\ar@{-}[dl]\\ & \bullet & }}
\end{align*}
According to Lemma \ref{lma:gauge}, we have 
$$
\frac{df}{dx}' x\otimes\frac{df}{dx}'' +\frac{df}{dy}'y\otimes
\frac{df}{dy}''+ 1\otimes f
= \frac{df}{dx}'\otimes x \frac{df}{dx}'' + \frac{df}{dy}'\otimes
y\frac{df}{dy}'' + f\otimes 1,
$$
and the same for $g$. Applying $-\circ \frac{d}{dx}$ (where $\circ$
means the right inner multiplication) and then the right (outer)
multiplication by $\frac{d}{dx}$, we obtain:
\begin{align*}
\raisebox{1.8pc}{\xymatrix@R=.75pc@C=.75pc{& \circ& \\
    \lck{\frac{df}{dx}'x}\ar@{-}[ur]\ar@{-}[dr] & &
    \frac{df}{dx}''\ar@{-}[ul]\ar@{-}[dl]\\ & \circ & }}
+\raisebox{1.8pc}{\xymatrix@R=.75pc@C=.75pc{& \circ& \\
    \lck{\frac{df}{dy}'y}\ar@{-}[ur]\ar@{-}[dr] & &
    \frac{df}{dy}''\ar@{-}[ul]\ar@{-}[dl]\\ & \circ & }}
+\raisebox{1.8pc}{\xymatrix@R=.75pc@C=.75pc{& \circ& \\
    \lck{1}\ar@{-}[ur]\ar@{-}[dr] & &
    f\ar@{-}[ul]\ar@{-}[dl]\\ & \circ & }} = \\
\raisebox{1.8pc}{\xymatrix@R=.75pc@C=.75pc{& \circ& \\
    \lck{\frac{df}{dx}'}\ar@{-}[ur]\ar@{-}[dr] & &
    x\frac{df}{dx}''\ar@{-}[ul]\ar@{-}[dl]\\ & \circ & }}
+\raisebox{1.8pc}{\xymatrix@R=.75pc@C=.75pc{& \circ& \\
    \lck{\frac{df}{dy}'}\ar@{-}[ur]\ar@{-}[dr] & &
    y\frac{df}{dy}''\ar@{-}[ul]\ar@{-}[dl]\\ & \circ & }}
+\raisebox{1.8pc}{\xymatrix@R=.75pc@C=.75pc{& \circ& \\
    \lck{f}\ar@{-}[ur]\ar@{-}[dr] & &
    1\,,\ar@{-}[ul]\ar@{-}[dl]\\ & \circ & }}
\end{align*}
whence
\begin{align*}
(A) = & 
\,-\raisebox{1.8pc}{\xymatrix@R=.75pc@C=.75pc{& \circ & \\ 
\lck{1}\ar@{-}[ur]\ar@{-}[dr] & & f\ar@{-}[ul]\ar@{-}[dl]\\ & \circ & }}
+ \raisebox{1.8pc}{\xymatrix@R=.75pc@C=.75pc{& \circ & \\ \lck{\frac{df}{dy}'}\ar@{-}[ur]\ar@{-}[dr] & & y\frac{df}{dy}''\,.\ar@{-}[ul]\ar@{-}[dl]\\ & \circ & }}
\end{align*}
A similar argument for $g$ (applying
$-\circ\frac{d}{dx}$ and then the right multiplication by
$\frac{d}{dy}$) yields
\begin{align*}
(B) = & 
\,-\raisebox{1.8pc}{\xymatrix@R=.75pc@C=.75pc{& \circ & \\ 
\lck{1}\ar@{-}[ur]\ar@{-}[dr] & & g\ar@{-}[ul]\ar@{-}[dl]\\ & \bullet & }}
+ \raisebox{1.8pc}{\xymatrix@R=.75pc@C=.75pc{& \circ & \\
    \lck{\frac{dg}{dy}'}\ar@{-}[ur]\ar@{-}[dr] & &
    y\frac{dg}{dy}''\ar@{-}[ul]\ar@{-}[dl]\\ & \bullet & }}
+ \raisebox{1.8pc}{\xymatrix@R=.75pc@C=.75pc{& \circ & \\ \lck{y\frac{df}{dx}''}\ar@{-}[ur]\ar@{-}[dr] & & \frac{df}{dx}'\,.\ar@{-}[ul]\ar@{-}[dl]\\ & \bullet & }}
\end{align*}
Adding the expressions obtained for $(A)$, $(B)$ and $(C)$, leads to the expression of $d_{P_1}^1(f\frac{d}{dx}+g\frac{d}{dy})$ stated in the lemma.
\end{proof}

Now we can consider the double Poisson cohomology space
$H^1_{P_1}(\Cxy)= \displaystyle \ker(d_{P_1}^1)/\hbox{Im}(d_{P_1}^0)$.
To do this, let $f\frac{d}{dx}+g\frac{d}{dy}\in\ker(d_{P_1}^1)$ be a
$1$-cocycle. According to Lemma \ref{imageofdp11}, this means
\begin{align}\label{A=0}
(A) =& 
\,-\raisebox{1.8pc}{\xymatrix@R=.75pc@C=.75pc{& \circ & \\ 
\lck{1}\ar@{-}[ur]\ar@{-}[dr] & & f\ar@{-}[ul]\ar@{-}[dl]\\ & \circ & }}
+ \raisebox{1.8pc}{\xymatrix@R=.75pc@C=.75pc{& \circ & \\ \lck{\frac{df}{dy}'}\ar@{-}[ur]\ar@{-}[dr] & & y\frac{df}{dy}''\ar@{-}[ul]\ar@{-}[dl]\\ & \circ & }}\,=0,
\end{align}
\begin{align}\label{B=0}
(B) =& 
\,-\raisebox{1.8pc}{\xymatrix@R=.75pc@C=.75pc{& \circ & \\ 
\lck{1}\ar@{-}[ur]\ar@{-}[dr] & & g\ar@{-}[ul]\ar@{-}[dl]\\ & \bullet & }}
+ \raisebox{1.8pc}{\xymatrix@R=.75pc@C=.75pc{& \circ & \\
    \lck{\frac{dg}{dy}'}\ar@{-}[ur]\ar@{-}[dr] & &
    y\frac{dg}{dy}''\ar@{-}[ul]\ar@{-}[dl]\\ & \bullet & }}
+ \raisebox{1.8pc}{\xymatrix@R=.75pc@C=.75pc{& \circ & \\
    \lck{y\frac{df}{dx}''}\ar@{-}[ur]\ar@{-}[dr] & &
    \frac{df}{dx}'\ar@{-}[ul]\ar@{-}[dl]\\ & \bullet & }} \, =0
\end{align}
and
\begin{align}\label{C=0}
(C)=& - \raisebox{1.8pc}{\xymatrix@R=.75pc@C=.75pc{& \bullet& \\ \lck{\frac{dg}{dx}'}\ar@{-}[ur]\ar@{-}[dr] & & y\frac{dg}{dx}''\ar@{-}[ul]\ar@{-}[dl]\\ & \bullet & }}\,=0.
\end{align}
Equation (\ref{A=0}) yields $f=y\tilde{f} + a$, with $a\in\CC$ and
$\tilde f \in\Cxy$. Indeed, write $f=y\tilde{f} + x \tilde{h}+
a$, with $\tilde f,\tilde h \in\Cxy$ and $a\in\CC$. Then, as 
$\raisebox{1.8pc}{\xymatrix@R=.75pc@C=.75pc{& \circ & \\ 
\lck{1}\ar@{-}[ur]\ar@{-}[dr] & & 1\ar@{-}[ul]\ar@{-}[dl]\\ & \circ &
}}=0$ (up to commutators), we have 
\begin{align*}
(A) =& 
-\raisebox{1.8pc}{\xymatrix@R=.75pc@C=.75pc{& \circ & \\ 
\lck{1}\ar@{-}[ur]\ar@{-}[dr] & & x\tilde{h}\ar@{-}[ul]\ar@{-}[dl]\\ & \circ & }}
+ \raisebox{1.8pc}{\xymatrix@R=.75pc@C=.75pc{& \circ & \\
    \lck{y\frac{d\tilde{f}}{dy}'}\ar@{-}[ur]\ar@{-}[dr] & &
    y\frac{d\tilde{f}}{dy}''\ar@{-}[ul]\ar@{-}[dl]\\ & \circ & }}
+ \raisebox{1.8pc}{\xymatrix@R=.75pc@C=.75pc{& \circ & \\
    \lck{x\frac{d\tilde{h}}{dy}'}\ar@{-}[ur]\ar@{-}[dr] & &
    y\frac{d\tilde{h}}{dy}''\ar@{-}[ul]\ar@{-}[dl]\\ & \circ & }}\,=0.
\end{align*}
For this equality to hold, the first term has to cancel itself, so that $\tilde{h}$ has to be
zero and $f=y\tilde{f} + a$.

A similar argument for Equation (\ref{B=0}) shows that $g=y\tilde{g}$, with $\tilde{g}\in\Cxy$ (but, in contrast with $f$, $g$ can not be a constant because
$\raisebox{1.8pc}{\xymatrix@R=.75pc@C=.75pc{& \circ & \\ 
\lck{1}\ar@{-}[ur]\ar@{-}[dr] & & 1\ar@{-}[ul]\ar@{-}[dl]\\ & \bullet &
}}\not=0$). Now, Equation (\ref{C=0}) becomes
\begin{align*}
(C)=& - \raisebox{1.8pc}{\xymatrix@R=.75pc@C=.75pc{& \bullet& \\ \lck{y\frac{d\tilde{g}}{dx}'}\ar@{-}[ur]\ar@{-}[dr] & & y\frac{d\tilde{g}}{dx}''\ar@{-}[ul]\ar@{-}[dl]\\ & \bullet & }}\,=0,
\end{align*}
so that 
$$
\frac{d\tilde{g}}{dx}= m'\otimes m'' + m''\otimes m',
$$
with $m',m''\in\Cxy$ (using Sweedler notation).
Using the NC-Euler Formula (Proposition \ref{nceuler}), we
can now write 
\begin{equation}\label{eqn:expn_tildeg}
\tilde{g} = \frac{1}{\deg_x(m'm'')+1}\left(m'xm'' + m''xm'\right) + p(y),
\end{equation}
where $p\in\CC\langle y\rangle$. Then, computing $\displaystyle\frac{d\tilde{g}}{dx}$
again, we get
\begin{eqnarray*}
\lefteqn{m'\otimes m'' + m''\otimes m'=}\\
&&\quad \frac{1}{\deg_x(m'm'')+1}\left(\frac{dm'}{dx}xm''
 +m'\otimes m''+ m'x\frac{dm''}{dx}\right.\\
&&\quad\qquad\qquad\left. +\, \frac{dm''}{dx}xm'+m''\otimes m' 
+ m''x\frac{dm'}{dx}\right),
\end{eqnarray*}
that is to say
\begin{eqnarray}\label{eqn:cond1}
\lefteqn{(\deg_x(m'm''))\left(m'\otimes m'' + m''\otimes
 m'\right)=}\nonumber\\
&& \quad\frac{dm'}{dx}xm'' + m'x\frac{dm''}{dx} + \frac{dm''}{dx}xm'+ m''x\frac{dm'}{dx}.
\end{eqnarray}
Now let $h=-\left(m''xm'x + m'xm''x\right)$ and $k=-p(y)x$ and let us compute
$\displaystyle -y\frac{dh}{dx}''\frac{dh}{dx}'$ and 
$\displaystyle -y\frac{dk}{dx}''\frac{dk}{dx}'$.
First, we have 
\begin{eqnarray*}
\frac{dh}{dx} &=& -\frac{dm''}{dx}x m'x - m'' \otimes m'x - m'' x
\frac{dm'}{dx}x - m''xm'\otimes 1\\
& -&\, \frac{dm'}{dx}x m''x - m' \otimes m''x - m' x \frac{dm''}{dx}x
- m'xm''\otimes 1,\\
\frac{dk}{dx}&=& -p\otimes 1,
\end{eqnarray*}
so that 
\begin{eqnarray*}
-y\frac{dh}{dx}''\frac{dh}{dx}' &=& 
2\, y\left(\frac{dm''}{dx}\right)''xm'x\left(\frac{dm''}{dx}\right)' + 2\,
ym'xm''\\
&+&\, 2\,
y\left(\frac{dm'}{dx}\right)''xm''x\left(\frac{dm'}{dx}\right)'+ 2\,
ym''xm',\\
-y\frac{dk}{dx}''\frac{dk}{dx}' &=& yp(y).
\end{eqnarray*}
But, applying the left outer multiplication by $x$, $-^{op}$ and $\mu$ to
Equation (\ref{eqn:cond1}), we obtain
\begin{eqnarray*}
\lefteqn{(\deg_x(m'm''))\left( m''xm' + m'x m''\right)=}\\
&&\qquad 2\left( \frac{dm'}{dx}\right)''xm''x\left(\frac{dm'}{dx}\right)' 
+ 2\left( \frac{dm''}{dx}\right)''xm'x\left(\frac{dm''}{dx}\right)'.
\end{eqnarray*}
This implies,
\begin{eqnarray*}
-y\frac{dh}{dx}''\frac{dh}{dx}' = 
 (\deg_x(m'm'')+2)\,\left( ym'xm''+ ym''xm'\right).
\end{eqnarray*}
In combination with (\ref{eqn:expn_tildeg}) we obtain
\begin{eqnarray*}
g &=& y\tilde{g} 
=
\frac{-1}{(\deg_x(m'm'')+1)(\deg_x(m'm'')+2)}\,y\,\frac{dh}{dx}''\frac{dh}{dx}' 
-y\frac{dk}{dx}''\frac{dk}{dx}'. 
\end{eqnarray*}
Now, we want to write $f$ in terms of $\displaystyle\frac{dh}{dy}$ 
and $\displaystyle\frac{dk}{dy}$. To do this, we
will use the Equation (\ref{B=0}). Using $f=y\tilde{f}+a$ and
$g=y\tilde{g}$, this equation can be written as follows:
\begin{align}\label{B=0'}
(B) = &\, 
 \raisebox{1.8pc}{\xymatrix@R=.75pc@C=.75pc{& \circ & \\
    \lck{y\frac{d\tilde{g}}{dy}'}\ar@{-}[ur]\ar@{-}[dr] & &
    y\frac{d\tilde{g}}{dy}''\ar@{-}[ul]\ar@{-}[dl]\\ & \bullet & }}
+ \raisebox{1.8pc}{\xymatrix@R=.75pc@C=.75pc{& \circ & \\
    \lck{y\frac{d\tilde{f}}{dx}''}\ar@{-}[ur]\ar@{-}[dr] & &
    y\frac{d\tilde{f}}{dx}'\ar@{-}[ul]\ar@{-}[dl]\\ & \bullet & }} \, =0.
\end{align}
This implies
$$
\frac{d\tilde{f}}{dx} = -\frac{d\tilde{g}}{dy}''\otimes\frac{d\tilde{g}}{dy}'.
$$
Using this expression and the NC-Euler Formula (Proposition
\ref{nceuler}), we get:
\begin{eqnarray*}
\tilde{f} &=& 
\frac{1}{\left(\deg_x\left(\frac{d\tilde{f}}{dx}'\frac{d\tilde{f}}{dx}''\right)+1\right)}\,
\frac{d\tilde{f}}{dx}'x\frac{d\tilde{f}}{dx}'' + l(y)\\
&=& \frac{-1}{\left(\deg_x\left(\frac{d\tilde{g}}{dy}'\frac{d\tilde{g}}{dy}''\right)+1\right)}\,
\frac{d\tilde{g}}{dy}''x\frac{d\tilde{g}}{dy}' + l(y),
\end{eqnarray*}
where $l\in\CC\langle y\rangle$. Now, the expression for $\tilde{g}$ obtained in 
(\ref{eqn:expn_tildeg}) yields
\begin{eqnarray*}
\tilde{f} &=& 
\frac{-2}{\left(\deg_x\left(\frac{dm'}{dy}m''\right)+2\right)
\left(\deg_x(m'm'')+1\right)}\,
\left(\left(\frac{dm'}{dy}\right)''xm''x\left(\frac{dm'}{dy}\right)'\right)\\
&-&
\frac{2}{\left(\deg_x\left(\frac{dm''}{dy}m'\right)+2\right)
\left(\deg_x(m'm'')+1\right)}\,
\left(\left(\frac{dm''}{dy}\right)''xm'x\left(\frac{dm''}{dy}\right)'\right)\\
&-&
\frac{1}{\left(\deg_x\left(\frac{dp}{dy}'\frac{dp}{dy}''\right)+1\right)}\,
\frac{dp}{dy}''x\frac{dp}{dy}' + l(y).
\end{eqnarray*}
Now, as $\deg_x\left(\displaystyle\frac{dm'}{dy}\right)=\deg_x(m')$ (unless
 $m'\in\CC\langle x\rangle$, that is, unless $\displaystyle\frac{dm'}{dy}=0$), which also holds for $m''$, we have exactly
\begin{eqnarray*}
\tilde{f} &=& 
\frac{-2}{\left(\deg_x(m'm'')+2\right)
\left(\deg_x(m'm'')+1\right)}\,
\left(\left(\frac{dm'}{dy}\right)''xm''x\left(\frac{dm'}{dy}\right)'\right)\\
&-&
\frac{2}{\left(\deg_x(m''m')+2\right)
\left(\deg_x(m'm'')+1\right)}\,
\left(\left(\frac{dm''}{dy}\right)''xm'x\left(\frac{dm''}{dy}\right)'\right)\\
&-&
\frac{dp}{dy}''x\frac{dp}{dy}' + l(y)\\
&=& \frac{1}{\left(\deg_x(m''m')+2\right)
\left(\deg_x(m'm'')+1\right)}\,\frac{dh}{dy}''\frac{dh}{dy}'\\
&+& \frac{dk}{dy}''\frac{dk}{dy}' + l(y).
\end{eqnarray*}
So that, if 
\begin{eqnarray*}
L &=& \frac{1}{\left(\deg_x(m''m')+2\right)
\left(\deg_x(m'm'')+1\right)}\, h\, +\, k \\
&=& \frac{-1}{\left(\deg_x(m''m')+2\right)
\left(\deg_x(m'm'')+1\right)}\,\left(m''xm'x + m'xm''x\right) -p(y)x,
\end{eqnarray*}
we have shown that
$$
g=-y\frac{dL}{dx}''\frac{dL}{dx}',\quad
f=y\frac{dL}{dy}''\frac{dL}{dy}' +y\,l(y) +a.
$$
Finally, as for every $n\in \NN^*$,
$y^n=y\frac{dq}{dy}''\frac{dq}{dy}'$, with 
$\displaystyle q=\frac{1}{n}\,y^n\in\CC\langle y\rangle $, the element
$y\,l(y)$ is of the form $\displaystyle y\frac{dQ}{dy}''\frac{dQ}{dy}'$, with
$Q\in\CC\langle y\rangle$ (and in particular
$\displaystyle y\frac{dQ}{dx}''\frac{dQ}{dx}'=0$) and
$$
f\frac{d}{dx} + g\frac{d}{dy} = d^0_{P_1}(L+Q) + a\frac{d}{dx}.
$$
As $\frac{d}{dx}\not\in \hbox{Im} d^0_{P_1}$, we have shown
\begin{prop}
The first double Poisson cohomology group of $\Cxy$, associated to the double
Poisson tensor 
$\displaystyle P_1=x\frac{d}{dx}\frac{d}{dx} + y\frac{d}{dx}\frac{d}{dy}$ is given
by
$$
H^1_{P_1}(\Cxy) \simeq \CC\,\frac{d}{dx}.
$$
\end{prop}
\begin{remark}
If we consider the double Poisson tensor 
$$
\tilde P_1 := P_{lin}^{B_2^2}  =  x\frac{\partial}{\partial
  x}\frac{\partial}{\partial y} + y\frac{\partial}{\partial
  y}\frac{\partial}{\partial y},
$$
we obtain in a similar fashion to the computations above for $P_1=P_{lin}^{B_2^1}$,
$$
H^0_{\tilde P_1} \simeq \CC \quad \hbox{and}\quad
H^1_{\tilde P_1} \simeq \CC\,\frac{d}{dy}.
$$
\end{remark}
Let us now consider the (classical) Poisson bracket on $\CC[x,y]$,
associated to $P_1$, that is $tr(P_1)=\pi_1 = y\,\frac{d}{dx}\wedge
\frac{d}{dy}$. According to \cite{Monnier} or \cite{RV}, or by a direct computation, we
have 
\begin{eqnarray*}
H^0_{\pi_1}(\CC[x,y]) &=& \CC, \qquad H^1_{\pi_1}(\CC[x,y]) = \CC
\frac{d}{dx},\\
H^i_{\pi_1}(\CC[x,y]) &=& 0, \qquad \hbox{for all } i\geq 2.
\end{eqnarray*}
So that the map $tr : H^i_{P_1}(\Cxy)\to H^i_{\pi_1}(\CC[x,y])$ is
bijective, for $i=0,1$.
\subsection{The nonlinear double Poisson tensor 
$P=x\frac{d}{dx}x\frac{d}{dy}$}

We conclude this section with the determination of the first two double Poisson cohomology groups of a nonlinear double Poisson bracket on the free algebra in two variables.
\begin{lemma}\label{thebracket}
The double bracket $\db{-,-}$ defined on $\Cxy$ as
$$\db{x,x} = \db{y,y} = 0 \mathrm{~and~} \db{x,y} = x\otimes x$$
is a double Poisson bracket.
\end{lemma}
\begin{proof}
First of all note that this bracket is defined by the double Poisson tensor $x\frac{d}{dx}x\frac{d}{dy}$. Representing $\frac{d}{dx}$ by $\circ$ and $\frac{d}{dy}$ by $\bullet$, this double Poisson bracket corresponds to the necklace $P$ depicted as
$$
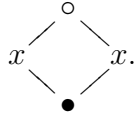
\xymatrix@R=.75pc@C=.75pc{& \circ\ar@{-}[dr] & \\ x\ar@{-}[ur]\ar@{-}[dr] & & x.\ar@{-}[dl] \\ & \bullet &}$$
The NC-Schouten bracket of $P$ with itself now becomes
$$
\raisebox{1.8pc}{\xymatrix@R=.75pc@C=.75pc{
& \circ \ar@{-}[r] & x \\
\lck{x}\ar@{-}[ur]\ar@{-}[dr] & & & \bullet\ar@{-}[ul]\ar@{-}[dl] \\
& \bullet\ar@{-}[r] & x
}}
+
\raisebox{1.8pc}{\xymatrix@R=.75pc@C=.75pc{
& \bullet \ar@{-}[r] & x \\
\lck{x}\ar@{-}[ur]\ar@{-}[dr] & & & \circ\ar@{-}[ul]\ar@{-}[dl] \\
& \bullet\ar@{-}[r] & x
}}
-
\raisebox{1.8pc}{\xymatrix@R=.75pc@C=.75pc{
& \circ \ar@{-}[r] & x \\
\lck{x}\ar@{-}[ur]\ar@{-}[dr] & & & \bullet\ar@{-}[ul]\ar@{-}[dl] \\
& \bullet\ar@{-}[r] & x
}}
-
\raisebox{1.8pc}{\xymatrix@R=.75pc@C=.75pc{
& \bullet \ar@{-}[r] & x \\
\lck{x}\ar@{-}[ur]\ar@{-}[dr] & & & \circ\ar@{-}[ul]\ar@{-}[dl] \\
& \bullet\ar@{-}[r] & x
}} = 0.
$$
\end{proof}

\begin{remark}
To stress the difference between double Poisson brackets and ordinary
Poisson brackets, note that $y\frac{d}{dx}y\frac{d}{dy}$ is also a
double Poisson tensor. However, taking higher degree monomials in $x$
or $y$ no longer yields double Poisson tensors. Whereas, of course,
for $\CC[x,y]$, any polynomial $\psi$ in $x$ and $y$ defines a Poisson
bracket $\psi\,\frac{d}{dx}\wedge\frac{d}{dy}$.
\end{remark}

For the remainder of this section, $P$ will be the double Poisson tensor $x\frac{d}{dx}x\frac{d}{dy}$.
First of all, observe that
\begin{prop}\label{imageofdp0}
For $f\in\Cxy$, we have
$$d_P(f) = 
\xymatrix@R=1.75pc{\lck{\circ}\ar@{-}@/^.25pc/[r] & x\frac{df}{dy}''\frac{df}{dy}'x\ar@{-}@/^.5pc/[l]}
-
\xymatrix@R=.75pc{\lck{\bullet}\ar@{-}@/^.25pc/[r] & x\frac{df}{dx}''\frac{df}{dx}'x\ar@{-}@/^.5pc/[l]}.
$$
Which means that
$$H_P^0(\Cxy) = \CC.$$
\end{prop}
\begin{proof}
The computation of $d_P(f)$ was already done in greater generality in Section \ref{intro}. To compute $H^0_P(\Cxy)$, note that 
$$\xymatrix@R=1.75pc{\lck{\circ}\ar@{-}@/^.25pc/[r] & x\frac{df}{dy}''\frac{df}{dy}'x\ar@{-}@/^.5pc/[l]}
-
\xymatrix@R=.75pc{\lck{\bullet}\ar@{-}@/^.25pc/[r] & x\frac{df}{dx}''\frac{df}{dx}'x\ar@{-}@/^.5pc/[l]}=0.
$$
implies
$$ \frac{df}{dx}''\frac{df}{dx}' = \frac{df}{dy}''\frac{df}{dy}' = 0.$$
But then
$$x\frac{df}{dx}''\frac{df}{dx}' + \frac{df}{dy}''\frac{df}{dy}'y = 0.$$
This means
$$\frac{df}{dx}'x\frac{df}{dx}'' + \frac{df}{dy}'y\frac{df}{dy}'' \in [\Cxy,\Cxy],$$
which by the NC-Euler formula (Proposition \ref{nceuler}) implies that
$$f\in \CC \oplus [\Cxy,\Cxy].$$
Now $H^0_P(\Cxy) = \ker(d_P^0)/[\Cxy,\Cxy]$, finishing the proof.
\end{proof}

Next, we can state that
\begin{lemma}\label{imageofdp1}
Let $f\frac{d}{dx}+g\frac{d}{dy}\in (T_{\Cxy}/[T_{\Cxy},T_{\Cxy}])_1$, then
\begin{align*}
d_P\left(f\frac{d}{dx}+g\frac{d}{dy}\right) = & 
- \raisebox{1.8pc}{\xymatrix@R=.75pc@C=.75pc{& \circ & \\ \lck{x}\ar@{-}[ur]\ar@{-}[dr] & & f\ar@{-}[ul]\ar@{-}[dl]\\ & \bullet & }}
- \raisebox{1.8pc}{\xymatrix@R=.75pc@C=.75pc{& \circ & \\ \lck{f}\ar@{-}[ur]\ar@{-}[dr] & & x\ar@{-}[ul]\ar@{-}[dl]\\ & \bullet & }}
+ \raisebox{1.8pc}{\xymatrix@R=.75pc@C=.75pc{& \circ & \\ \lck{\frac{df}{dy}'x}\ar@{-}[ur]\ar@{-}[dr] & & x\frac{df}{dy}''\ar@{-}[ul]\ar@{-}[dl]\\ & \circ & }}  \\
&
+ \raisebox{1.8pc}{\xymatrix@R=.75pc@C=.75pc{& \circ & \\ \lck{x\frac{df}{dx}''}\ar@{-}[ur]\ar@{-}[dr] & & \frac{df}{dx}'x\ar@{-}[ul]\ar@{-}[dl]\\ & \bullet & }}
+ \raisebox{1.8pc}{\xymatrix@R=.75pc@C=.75pc{& \circ & \\ \lck{\frac{dg}{dy}'x}\ar@{-}[ur]\ar@{-}[dr] & & x\frac{dg}{dy}''\ar@{-}[ul]\ar@{-}[dl]\\ & \bullet & }}
- \raisebox{1.8pc}{\xymatrix@R=.75pc@C=.75pc{& \bullet& \\ \lck{\frac{dg}{dx}'x}\ar@{-}[ur]\ar@{-}[dr] & & x\frac{dg}{dx}''\ar@{-}[ul]\ar@{-}[dl]\\ & \bullet & }}\,.
\end{align*}
We will denote this expression by $(*)$.
\end{lemma}
\begin{proof}
Straightforward.
\end{proof}
Using this lemma, we can determine the kernel of $d_P^1$. So assume
now that $f\frac{d}{dx}+g\frac{d}{dy}\in\ker(d_P^1)$. First of all
note that if $d_P\left(f\frac{d}{dx}+g\frac{d}{dy}\right) = 0$, the
third term in the expression $(*)$ for
$d_P\left(f\frac{d}{dx}+g\frac{d}{dy}\right)$ in Lemma
\ref{imageofdp1} has to cancel itself and the sixth term in this
expression has to cancel itself. This implies (using the Sweedler notations) 
$$\frac{df}{dy} = xm_f'\otimes m_f''x + xm_f''\otimes m_f'x + n_f' \otimes n_f''x + xn_f''\otimes n_f' + c_f 1\otimes 1$$
with $m_f',m_f'',n_f''\in \Cxy$, $c_f, n_f'\in\CC$
and
$$\frac{dg}{dx} = xm_g'\otimes m_g''x + xm_g''\otimes m_g'x + n_g' \otimes n_g''x + xn_g''\otimes n_g' + c_g 1\otimes 1$$
with $m_g',m_g'',n_g''\in \Cxy$, $c_g, n_g'\in\CC$.

Using the NC-Euler formula (Proposition \ref{nceuler}), this implies
\begin{eqnarray*}
f &=& \frac{1}{\deg_y(m_f'm_f'')+1}\,x(m_f' y m_f'' + m_f'' y m_f')x \\
&&+ \frac{1}{\deg_y(n_f'')+1}\,(n_f'  y n_f''x + xn_f'' y n_f') + p(x) +
c_fy
\end{eqnarray*}
and
\begin{eqnarray*}
g &=& \frac{1}{\deg_x(m_g'm_g'')+3}\,x(m_g' x m_g'' + m_g'' x m_g')x \\
&&+\frac{2n_g'}{\deg_x(n_g'')+2} \,x n_g''x + q(y) + c_gx.
\end{eqnarray*}
Now note that for $c_fy$, the first two terms of $(*)$ yield
$$
- c_f\,\raisebox{1.8pc}{\xymatrix@R=.75pc@C=.75pc{& \circ & \\ \lck{x}\ar@{-}[ur]\ar@{-}[dr] & & y\ar@{-}[ul]\ar@{-}[dl]\\ & \bullet & }}
- c_f\,\raisebox{1.8pc}{\xymatrix@R=.75pc@C=.75pc{& \circ & \\ \lck{y}\ar@{-}[ur]\ar@{-}[dr] & & x\ar@{-}[ul]\ar@{-}[dl]\\ & \bullet & }}.
$$
Because of the degree in $x$ of the remaining terms, these terms cannot vanish unless $c_f = 0$.

Using this last remark and the expression for $f$ above to compute $\frac{df}{dy}$ again, we obtain
\begin{align*}
\frac{df}{dy} = & \frac{1}{\deg_y(m_f'm_f'')+1}\,x\left(\frac{dm_f'}{dy} y m_f'' + m_f'\otimes m_f'' + m_f'y\frac{dm_f''}{dy} + \frac{dm_f''}{dy} y m_f'\right. \\
&\left. + m_f''\otimes m_f' + m_f''y\frac{dm_f'}{dy}\right)x \\
& + \frac{1}{\deg_y(n_f'')+1}(n_f'  \otimes n_f''x + n_f'y\frac{dn_f''}{dy}x+ x\frac{dn_f''}{dy} y n_f' + xn_f''\otimes n_f').
\end{align*}
This expression should be equal to the first expression found for $\frac{df}{dy}$. That is,
\begin{align*}
& xm_f'\otimes m_f''x + xm_f''\otimes m_f'x + n_f' \otimes n_f''x + xn_f''\otimes n_f'  \\
= & \frac{1}{\deg_y(m_f'm_f'')+1}\,x\left(\frac{dm_f'}{dy} y m_f'' + m_f'\otimes m_f'' + m_f'y\frac{dm_f''}{dy} + \frac{dm_f''}{dy} y m_f' \right.\\
&\left. +\, m_f''\otimes m_f' + m_f''y\frac{dm_f'}{dy}\right)x \\
& + \frac{1}{\deg_y(n_f'')+1}(n_f'  \otimes n_f''x + n_f'y\frac{dn_f''}{dy}x+ x\frac{dn_f''}{dy} y n_f' + xn_f''\otimes n_f').
\end{align*}
whence, by comparing elements of the form $x\dots x$ we obtain
$$
\deg_y(m_f'm_f'') (m_f'\otimes m_f'' + m_f''\otimes m_f') = 
\frac{dm_f'}{dy} y m_f'' + m_f'y\frac{dm_f''}{dy} + \frac{dm_f''}{dy} y m_f' + m_f''y\frac{dm_f'}{dy}$$
and $n_f''\in\CC[x]$.

Letting $y$ act on the equality obtained in the previous paragraph by the left outer action, we obtain
\begin{eqnarray*}
\deg_y(m_f'm_f'') (ym_f'\otimes m_f'' + ym_f''\otimes m_f') &=& 
y\frac{dm_f'}{dy} y m_f'' + ym_f'y\frac{dm_f''}{dy}\\
&& +\, y\frac{dm_f''}{dy} y m_f' + ym_f''y\frac{dm_f'}{dy}.
\end{eqnarray*}
Which yields, using $-^{op}$ and $\mu$, the equality
\begin{eqnarray*}
\lefteqn{\deg_y(m_f'm_f'') (m_f''ym_f' + m_f'ym_f'') = }\\
&&2\left(   \left(\frac{dm_f'}{dy}\right)''ym_f''y\left(\frac{dm_f'}{dy}\right)' 
+
\left(\frac{dm_f''}{dy}\right)''ym_f'y\left(\frac{dm_f''}{dy}\right)'\right).
\end{eqnarray*}
Now let $h = \frac{2}{(\deg_y(m_f'm_f'')+2)\,(\deg_y(m_f'm_f'')+1)} m_f''ym_f'y$, then
\begin{align*}
x\frac{dh}{dy}''\frac{dh}{dy}'x = &  \frac{1}{\deg_y(m_f'm_f'')+1}(xm_f'ym_f''x  + xm_f''ym_f'x) 
\end{align*}
and
\begin{eqnarray*}
x\frac{dh}{dx}''\frac{dh}{dx}'x &=& \frac{2}{(\deg_y(m_f'm_f'')+2)(\deg_y(m_f'm_f'')+1)}x\left( \left(\frac{dm_f''}{dx}\right)''ym_f'y \left(\frac{dm_f''}{dx}\right)' \right. \\
& & \left. + \left(\frac{dm_f'}{dx}\right)''ym_f''y\left(\frac{dm_f'}{dx}\right)'\right)x.
\end{eqnarray*}
So, writing
$$f_1 := f - x\frac{dh}{dy}''\frac{dh}{dy}'x :=  y p_1(x) + p_1(x) y + p(x)$$
with $p_1 := \sum_{i=1}^{n}a_ix^i$ and $p=\sum_{i=0}^mb_ix^i$ and
$$g_1 := g + x\frac{dh}{dx}''\frac{dh}{dx}'x,$$
we again obtain an element $f_1\frac{d}{dx}+ g_1\frac{d}{dy}$ in $\ker(d_P^1)$ by Proposition \ref{imageofdp0}. Observe moreover that $x^i$ for $i\geq 2$ can be written as $x\frac{dh_i}{dy}''\frac{dh_i}{dy}'x$ with $h_i = x^{i-2}y$, so we may assume (modifying $f$ and $g$ by the image of a suitable $h$)  $p(x) = b_1 x + b_0$. Now note that $b_0$ has to be equal to zero as only the first two terms of $(*)$ contain $b_0$  and these do not cancel each other. That is, we may assume $p(x) = b_1x$.

The image under $d_P^1$ of this element becomes
\begin{align*}
d_P\left(f_1\frac{d}{dx} + g_1\frac{d}{dy}\right) = & 
- \sum_{i=1}^na_i\, \raisebox{1.8pc}{\xymatrix@R=.75pc@C=.75pc{& \circ & \\ \lck{yx^i}\ar@{-}[ur]\ar@{-}[dr] & & x\ar@{-}[ul]\ar@{-}[dl]\\ & \bullet & }}
- \sum_{i=1}^na_i\,\raisebox{1.8pc}{\xymatrix@R=.75pc@C=.75pc{& \circ & \\ \lck{x}\ar@{-}[ur]\ar@{-}[dr] & & x^iy\ar@{-}[ul]\ar@{-}[dl]\\ & \bullet & }} \\
&
+\sum_{i=2}^n a_i\sum_{j=1}^{i-1} \raisebox{1.8pc}{\xymatrix@R=.75pc@C=.75pc{& \circ & \\ \lck{x^{i-j+1}}\ar@{-}[ur]\ar@{-}[dr] & & yx^j\ar@{-}[ul]\ar@{-}[dl]\\ & \bullet & }}
+\sum_{i=2}^n a_i\sum_{j=2}^i \raisebox{1.8pc}{\xymatrix@R=.75pc@C=.75pc{& \circ & \\ \lck{x^{i-j+1}y}\ar@{-}[ur]\ar@{-}[dr] & & x^j\ar@{-}[ul]\ar@{-}[dl]\\ & \bullet & }} \\
&
- b_1\raisebox{1.8pc}{\xymatrix@R=.75pc@C=.75pc{& \circ & \\ \lck{x}\ar@{-}[ur]\ar@{-}[dr] & & x\ar@{-}[ul]\ar@{-}[dl]\\ & \bullet & }} 
 + \raisebox{1.8pc}{\xymatrix@R=.75pc@C=.75pc{& \circ & \\ \lck{\frac{dg_1}{dy}'x}\ar@{-}[ur]\ar@{-}[dr] & & x\frac{dg_1}{dy}''\ar@{-}[ul]\ar@{-}[dl]\\ & \bullet & }}
- \raisebox{1.8pc}{\xymatrix@R=.75pc@C=.75pc{& \bullet& \\ \lck{\frac{dg_1}{dx}'x}\ar@{-}[ur]\ar@{-}[dr] & & x\frac{dg_1}{dx}''\ar@{-}[ul]\ar@{-}[dl]\\ & \bullet & }}
\end{align*}
Note that we cancelled two terms
$$
- \sum_{i=1}^n a_i\, \raisebox{1.8pc}{\xymatrix@R=.75pc@C=.75pc{& \circ & \\ \lck{x^iy}\ar@{-}[ur]\ar@{-}[dr] & & x\ar@{-}[ul]\ar@{-}[dl]\\ & \bullet & }}
- \sum_{i=1}^na_i\,\raisebox{1.8pc}{\xymatrix@R=.75pc@C=.75pc{& \circ & \\ \lck{x}\ar@{-}[ur]\ar@{-}[dr] & & yx^i\ar@{-}[ul]\ar@{-}[dl]\\ & \bullet & }}
$$
against the similar terms obtained in the second row of $(*)$ for $j = i$ in the first sum and $j=1$ in the second sum.

Now observe that for $n\geq 2$, the terms in the second row of the expression above cannot be eliminated by any other term because of the location of the $y$ factor, whence $a_i = 0$ for $i\geq 2$. That is, $f_1 = a(xy + yx) + p(x)$. Moreover, if $a\neq 0$, the expression
$$
-a\left(\raisebox{1.8pc}{\xymatrix@R=.75pc@C=.75pc{& \circ & \\ \lck{yx}\ar@{-}[ur]\ar@{-}[dr] & & x\ar@{-}[ul]\ar@{-}[dl]\\ & \bullet & }}
+ \raisebox{1.8pc}{\xymatrix@R=.75pc@C=.75pc{& \circ & \\ \lck{x}\ar@{-}[ur]\ar@{-}[dr] & & xy\ar@{-}[ul]\ar@{-}[dl]\\ & \bullet & }}\right)
$$
can only be cancelled if $g_1 = g_2 + ay^2$. That is, the image becomes
\begin{align*}
d_P\left(f_1\frac{d}{dx}+g_1\frac{d}{dy}\right) = & 
- b_1 \raisebox{1.8pc}{\xymatrix@R=.75pc@C=.75pc{& \circ & \\ \lck{x}\ar@{-}[ur]\ar@{-}[dr] & & x\ar@{-}[ul]\ar@{-}[dl]\\ & \bullet & }}
 + \raisebox{1.8pc}{\xymatrix@R=.75pc@C=.75pc{& \circ & \\ \lck{\frac{dg_2}{dy}'x}\ar@{-}[ur]\ar@{-}[dr] & & x\frac{dg_2}{dy}''\ar@{-}[ul]\ar@{-}[dl]\\ & \bullet & }}
- \raisebox{1.8pc}{\xymatrix@R=.75pc@C=.75pc{& \bullet& \\ \lck{\frac{dg_2}{dx}'x}\ar@{-}[ur]\ar@{-}[dr] & & x\frac{dg_2}{dx}''\ar@{-}[ul]\ar@{-}[dl]\\ & \bullet & }}
\end{align*}
Now if $b_1 \neq 0$, this expression can only be zero if $g_2 = g_3 +b_1y$, and we get
\begin{align*}
d_P\left(f_1\frac{d}{dx}+g_1\frac{d}{dy}\right) = & 
\, \raisebox{1.8pc}{\xymatrix@R=.75pc@C=.75pc{& \circ & \\ \lck{\frac{dg_3}{dy}'x}\ar@{-}[ur]\ar@{-}[dr] & & x\frac{dg_3}{dy}''\ar@{-}[ul]\ar@{-}[dl]\\ & \bullet & }}
- \raisebox{1.8pc}{\xymatrix@R=.75pc@C=.75pc{& \bullet& \\ \lck{\frac{dg_3}{dx}'x}\ar@{-}[ur]\ar@{-}[dr] & & x\frac{dg_3}{dx}''\ar@{-}[ul]\ar@{-}[dl]\\ & \bullet & }}
\end{align*}
However, the first term in this expression can only be zero if $\frac{dg_3}{dy} = 0$, so $g_3 \in \CC[x]$. Finally, observe that in $g_3$, we can cancel all monomials $x^i$ with $i\geq 2$ using $h = x^{i-1}$ (which does not modify $f$ in any way).

But then we have shown that
\begin{theorem}\label{descriptionofH1}
For $P$ as above, we have that
$$
H^1_P(\Cxy) \simeq \{\left(a(xy+yx) + bx\right)\frac{d}{dx}
+ \left(ay^2 + by + cx + d\right)\frac{d}{dy}\mid a,b,c,d\in\CC\},
$$
so in particular $\dim H^1_P(\Cxy) = 4$
\end{theorem}

\bigskip 

Let us now consider the Poisson bracket that corresponds to the double
Poisson tensor $P=x\frac{d}{dx}x\frac{d}{dy}$. We then obtain the
Poisson algebra $(\CC[x,y],\pi)$, where 
$\pi=tr(P) = x^2\frac{d}{dx}\wedge\frac{d}{dy}$.

According to \cite{RV}, as the polynomial $x^2$ is not square-free,
the first Poisson cohomology space $H^1_\pi(\CC[x,y])$
is infinite dimensional, so that 
\begin{cor}
The map $H^1_P(\Cxy)\rightarrow H^1_{tr(P)}(\CC[x,y])$ is not onto.
\end{cor}

Let us give an explicit basis for this vector space
$H^1_\pi(\CC[x,y])$, in order to make explicit this map.

First of all, we recall that the Poisson coboundary operator is given
by:
$d^k_\pi = \{ \pi, - \} : 
\bigwedge^k \Der(\CC[x,y])\to\bigwedge^{k+1} \Der(\CC[x,y])$, where
$\{-,-\}$ denotes the (classical) Schouten-Nijenhuis bracket and
$\Der(\CC[x,y])$ denotes the $\CC[x,y]$-module of the derivations of
the commutative algebra $\CC[x,y]$. 

We have , for $f,g,h\in\CC[x,y]$,
$$
\begin{array}{rcl}
d^0_\pi(h)&=&\displaystyle
x^2\left(-\frac{dh}{dy}\frac{d}{dx}+\frac{dh}{dx}\frac{d}{dy}\right)\\
d^1_\pi(f\frac{d}{dx}+g\frac{d}{dy}) 
&=&\displaystyle\left( x^2\left(\frac{df}{dx}+\frac{dg}{dy}\right)-2xf\right)\frac{d}{dx}\wedge\frac{d}{dy}.
\end{array}
$$
So that,
$$
\begin{array}{rcl}
H^0_\pi(\CC[x,y])\simeq \displaystyle\left\lbrace h\in\CC[x,y]\mid 
\frac{dh}{dy}=\frac{dh}{dx}=0\right\rbrace
\simeq \CC,\\
H^1_\pi(\CC[x,y])\simeq \frac{\displaystyle \Bigl\lbrace (f,g)\in \CC[x,y]^2
            \mid x^2\left(\frac{df}{dx}+\frac{dg}{dy}\right)-2xf =0\Bigr\rbrace}
                     {\displaystyle\Bigl\lbrace
     x^2\left(-\frac{dh}{dy},\frac{dh}{dx}\right) \mid h\in\CC[x,y]\Bigr\rbrace}.
\end{array}
$$
It is
clear that the coboundary operator $d^k_\pi$ is an homogeneous
operator, for example, if $f$ and $g$ are homogeneous polynomial of same degree
$d\in\NN$, then $d^1_\pi(f\frac{d}{dx}+g\frac{d}{dy})$ is given by
an homogeneous polynomial of degree $d+1$, in factor of
$\frac{d}{dx}\wedge\frac{d}{dy}$. This implies that we can work
``degree by degree'' and consider only homogeneous polynomials. We
recall the (commutative) Euler formula, for an homogeneous polynomial
$q\in\CC[x,y]$:
\begin{equation}\label{euler}
x\frac{dq}{dx} +y\frac{dq}{dy}= \deg(q)\,q.
\end{equation}

Let us consider $(f,g)\in\CC[x,y]^2$, two homogeneous polynomials of
same degree $d\in\NN$, satisfying the $1$-cocycle
condition $\displaystyle x^2\left(\frac{df}{dx}+\frac{dg}{dy}\right)=2xf$,
equivalent to $\displaystyle
x\left(\frac{df}{dx}+\frac{dg}{dy}\right)=2f$.
We divide $f$ and $g$ by $x^2$ and obtain:
$$
f=x^2 f_1 + x f_2 +f_3,\quad g=x^2 g_1 + x g_2 +g_3,
$$ 
with $f_1,g_1\in\CC[x,y]$ and $f_2,f_3,g_2,g_3\in\CC[y]$, homogeneous
polynomials.
Then the $1$-cocycle condition becomes:
$$
x\left(x^2 \frac{df_1}{dx} + f_2 + x^2 \frac{dg_1}{dy} +
x\frac{dg_2}{dy} + \frac{dg_3}{dy}\right)= 2x f_2 +2f_3,
$$
that leads to $f_3=0$ (because $f_3\in\CC[y]$) and 
$$
x^2 \frac{df_1}{dx} + f_2 + x^2 \frac{dg_1}{dy} +
x\frac{dg_2}{dy} + \frac{dg_3}{dy}= 2 f_2.
$$
We then have also $f_2 = \displaystyle\frac{dg_3}{dy}$ and 
$\displaystyle x \frac{df_1}{dx} + x \frac{dg_1}{dy} + \frac{dg_2}{dy}= 0$.
This equation then implies that $\displaystyle\frac{dg_2}{dy}=0$, i.e., $g_2\in\CC$
and also $\displaystyle\frac{df_1}{dx}=- \frac{dg_1}{dy}$.
Suppose now that $d\geq 2$ and let us consider the polynomial $h:=yf_1-xg_1$. 
We have, using the Euler formula (\ref{euler}) and the last equation above, 
\begin{eqnarray*}
\frac{dh}{dx} &=& y\frac{df_1}{dx} - g_1 -x\frac{dg_1}{dx}
 = -y\frac{dg_1}{dy} - g_1 -x\frac{dg_1}{dx}
= -(d-1) g_1,\\
\frac{dh}{dy} &=& f_1 + y\frac{df_1}{dy}  -x\frac{dg_1}{dy}
 = f_1 + y\frac{df_1}{dy}  +x\frac{df_1}{dx}
= (d-1) f_1.
\end{eqnarray*}
We have obtained that, if $d\geq 2$, then 
$x^2(f_1\frac{d}{dx} + g_1\frac{d}{dy})= d^0_\pi(-h)$. Moreover,
$g_3$ is an homogeneous polynomial of degree $d$, in $\CC[y]$, so that
$g_3=c_3y^d$, with $c_3\in\CC$. We have also seen that $f_3=0$,
$\displaystyle f_2=\frac{dg_3}{dy}=c_3\, d\, y^{d-1}$ and $g_2=c_2\in\CC$. 

It remains to consider the cases where $d=1$ and $d=0$. First, if
$d=0$, i.e., $f,g\in\CC$, then the $1$-cocycle condition is equivalent to
$f=0$, second, if $d=1$, we have $(f,g)=(ax+by,cx+dy)$, with
$a,b,c,d\in\CC$ and the $1$-cocycle condition says that $a=d$ and
$b=0$. We finally have obtain the following
\begin{prop}
The first Poisson cohomology space associated to the Poisson algebra 
$(\CC[x,y],\pi=x^2\frac{d}{dx}\wedge\frac{d}{dy})$ is given by:
\begin{equation*}
H^1_\pi(\CC[x,y])\simeq \bigoplus_{k\in \NN}\CC\left(k\,
y^{k-1}x,y^k\right) \oplus \CC(0,x). 
\end{equation*}
\end{prop}
The image of the double Poisson cohomology under the canonical trace map in the classical cohomology now becomes
\begin{cor}
For the double Poisson tensor $P = x\frac{\partial}{\partial x}x\frac{\partial}{\partial y}$ we have
$$H_{tr(P)}^1(\CC[x,y]) = \bigoplus_{k\geq 3}\CC\left(k
y^{k-1}x,y^k\right)\oplus tr(H_P^1(\Cxy)),$$
that is
$$tr(H_P^1(\Cxy)) = \CC\left(2yx,y^2\right)\oplus\CC(x,y)\oplus\CC(0,x)\oplus\CC(0,1).$$
\end{cor}


\begin{thebibliography}{10}

\bibitem{RafLieven}
Raf Bocklandt and Lieven Le~Bruyn.
\newblock Necklace {L}ie algebras and noncommutative symplectic geometry.
\newblock {\em Math. Z.}, 240(1):141--167, 2002.

\bibitem{CBEG}
Bill Crawley-Boevey, Pavel Etingov, and Victor Ginzburg.
\newblock Noncommutative geometry and quiver algebras.
\newblock math.AG/0502301, 2005.

\bibitem{Bill}
William Crawley-Boevey.
\newblock Geometry of the moment map for representations of quivers.
\newblock {\em Compositio Math.}, 126(3):257--293, 2001.

\bibitem{CQ}
Joachim Cuntz and Daniel Quillen.
\newblock Algebra extensions and nonsingularity.
\newblock {\em J. Amer. Math. Soc.}, 8(2):251--289, 1995.

\bibitem{Fuks}
D.~B. Fuks.
\newblock {\em Cohomology of infinite-dimensional {L}ie algebras}.
\newblock Contemporary Soviet Mathematics. Consultants Bureau, New York, 1986.
\newblock Translated from the Russian by A. B. Sosinski\u\i.

\bibitem{gab}
Peter Gabriel.
\newblock Finite representation type is open.
\newblock In {\em Proceedings of the International Conference on
  Representations of Algebras (Carleton Univ., Ottawa, Ont., 1974), Paper No.
  10}, pages 23 pp. Carleton Math. Lecture Notes, No. 9, Ottawa, Ont., 1974.
  Carleton Univ.

\bibitem{lazaroiu}
Calin~Iuliu Lazaroiu.
\newblock On the non-commutative geometry of topological {D}-branes.
\newblock {\em J. High Energy Phys.}, (11):032, 57 pp. (electronic), 2005.

\bibitem{Lich}
Andr{\'e} Lichnerowicz.
\newblock Les vari\'et\'es de {P}oisson et leurs alg\`ebres de {L}ie
  associ\'ees.
\newblock {\em J. Differential Geometry}, 12(2):253--300, 1977.

\bibitem{Monnier}
Philippe Monnier.
\newblock Poisson cohomology in dimension two.
\newblock {\em Israel J. Math.}, 129:189--207, 2002.

\bibitem{RV}
Claude Roger and Pol Vanhaecke.
\newblock Poisson cohomology of the affine plane.
\newblock {\em J. Algebra}, 251(1):448--460, 2002.

\bibitem{DPSSA}
Geert Van~de Weyer.
\newblock Double {P}oisson structures on finite-dimensional semi-simple
  algebras.
\newblock math.AG/0603533. To appear in Algebras and Representation Theory,
  2006.

\bibitem{MichelDPA}
Michel Van~den Bergh.
\newblock Double {P}oisson algebras.
\newblock math.AG/0410528, 2006 (updated).

\end{thebibliography}
\end{document}